\documentclass[VANCOUVER,STIX1COL]{WileyNJD-v2}

\articletype{SPECIAL ISSUE PAPER}%

\received{26 April 2016}
\revised{6 June 2016}
\accepted{6 June 2016}

\usepackage{diagbox}
\usepackage{float}
\usepackage{tikz}
\usepackage{slashbox}
\usepackage{bm}
\newcommand{\tikzcircle}[2][black,fill=black]{\tikz[baseline=-0.5ex]\draw[#1,radius=#2] (0,0) circle ;}%
\newcommand{\revise}[1]{{\color{black}#1}}

\usetikzlibrary{shapes,patterns,snakes}

\begin{document}

\title{A local Fourier analysis of additive Vanka relaxation for the Stokes equations}

\author[1]{Patrick E. Farrell}

\author[2]{Yunhui He*}

\author[2]{Scott P. MacLachlan}

\authormark{Farrell, He, and MacLachlan}

\address[1]{\orgdiv{Mathematical Institute}, \orgname{University of Oxford}, \orgaddress{\country{UK}}}

\address[2]{\orgdiv{Department of Mathematics and Statistics}, \orgname{Memorial University of Newfoundland}, \orgaddress{\state{St. John's, NL}, \country{Canada}}}

\corres{*Yunhui He, Department of Mathematics and Statistics, Memorial
  University of Newfoundland, St. John's, NL, A1C 5S7, Canada. \email{yunhui.he@mun.ca}}

\abstract[Summary]{Multigrid methods are popular solution algorithms for many discretized  PDEs, either as standalone iterative solvers or as preconditioners, due to their high efficiency. However, the choice and optimization of multigrid components such as relaxation schemes and grid-transfer operators is crucial to the design of optimally efficient algorithms. It is well--known that local Fourier analysis (LFA) is a useful tool to predict and analyze the performance of these components. In this paper, we develop a local Fourier analysis of monolithic multigrid methods based on additive Vanka relaxation schemes for mixed finite-element discretizations of  the Stokes equations. The analysis offers insight into the choice of  ``patches'' for the  Vanka relaxation, revealing that smaller patches offer more effective convergence per floating point operation. Parameters that minimize the two-grid convergence factor are proposed and numerical experiments are presented to validate the LFA predictions.}


\jnlcitation{\cname{%
\author{P.~E.~Farrell},
\author{Y.~He}, and
\author{S.~MacLachlan}} (\cyear{2019}),
\ctitle{A local Fourier analysis for additive Vanka relaxation for the Stokes equations}, \cjournal{Numer. Linear Alg. Appl.}, \cvol{2019;00:0--0}.}

\maketitle

%
%
\section{Introduction}\label{sec:intro}

Saddle-point problems are ubiquitous in applied mathematics.\cite{benzi2005numerical} Their importance motivates the development of effective parallel solvers. Block preconditioners and monolithic multigrid methods are established approaches for solving linear (or linearized) saddle-point problems. Block preconditioning is highly effective when the Schur complement of the system is well understood; for the Stokes equations, the Schur complement is spectrally equivalent to a weighted mass matrix, forming the basis for efficient solvers that use multigrid for the viscous term.\cite{wathen1993fast,silvester1994fast}  Monolithic  methods that apply multigrid to the entire system can also offer superb efficiency if an appropriate relaxation can be devised. For example, Adler et al.~\cite{JAdler_etal_2015b}  proposed a monolithic multigrid method with Braess-Sarazin relaxation for the Stokes equations that provided the fastest time to solution when compared with several block preconditioners and other monolithic multigrid methods.  While many block preconditioners have been successfully employed in (massively) parallel computing environments, \cite{TGeenen_etal_2009a, JRudi_etal_2015a} the same cannot be said for monolithic multigrid, whose parallelization was largely absent from the literature until recently. \cite{VJohn_LTobiska_2000a, BGmeiner_etal_2015a, BGmeiner_etal_2016a}

Common approaches to multigrid relaxation for coupled systems are  distributive relaxation \cite{niestegge1990analysis,oosterlee2006multigrid} (which  relies on continuum commutativity that may not hold at the discrete level), Braess-Sarazin relaxation, \cite{DBraess_RSarazin_1997a,adler2016monolithic} Uzawa-type relaxation,\cite{MR3217219,MR3592469} and Vanka relaxation. \cite{SPVanka_1986a,MR3848566} Based on distributive relaxation, Wang and Chen \cite{wang2013multigrid}  developed a least squares commutator distributive Gauss-Seidel  relaxation  for the Stokes equations. Furthermore, this technique has been extended to the Oseen problem by Chen et al. \cite{chen2015multigrid}  Braess-Sarazin relaxation is known to be highly efficient, and has been applied to nematic liquid crystals,\cite{MR3439774} magnetohydrodynamics, \cite{adler2016monolithic} and other coupled systems. Considering parallel computation, recently, He and MacLachlan  presented a local Fourier analysis (LFA) for both distributive weighted Jacobi and Braess-Sarazin relaxations for the Stokes equations discretized by the Marker-and-Cell finite-difference scheme \cite{he2018local} and  by mixed finite-element methods, \cite{HMFEMStokes} showing the power of LFA for designing efficient algorithms.

Vanka-type relaxation has been used in many contexts, such as for the Navier--Stokes equations, \cite{SPVanka_1986a,john2002higher} and extended to Vanka--like
  schemes for other problems or to improve performance. \cite{JAdler_etal_2015b, adler2016monolithic,VJohn_LTobiska_2000a, VJohn_GMatthies_2001a, SManservisi_2006a,schoberl2003schwarz, MR4024766, KahlKintscher_2018a, Claus_2019} However, Vanka relaxation is typically considered in its multiplicative variant. This seems overly constraining, particularly in
  comparison to Braess-Sarazin relaxation, which can naturally be done in additive form.\cite{he2018local,HMFEMStokes}  While multiplicative Vanka relaxation is very efficient, \revise{it has a higher cost per iteration than the additive variant and its parallelization is more involved}. We therefore consider additive variants of Vanka-type relaxation in this work.

There are two challenging choices to be made for Vanka relaxation, which are observed to be more critical in the additive setting.  First,  many  choices are possible for the underlying patches within the overlapping Schwarz framework. While we would naturally choose small patches for efficiency or large patches for effectiveness, no general results are known. Secondly, relaxation weights play an important role in ensuring best possible performance of the multigrid algorithm, particularly for additive methods.  Thus, there is a need for  analysis to inform the algorithmic  choices, and LFA seems well-suited. LFA has already been applied to Vanka relaxation in the multiplicative\cite{sivaloganathan1991use,molenaar1991two,SPMacLachlan_CWOosterlee_2011a,rodrigo2016local,MR4024766} and multicoloured\cite{boonen2008local,KahlKintscher_2018a} contexts; here, in contrast, we aim to develop LFA for additive schemes and use it to drive  parameter choices in practical experiments for the Stokes equations. To our knowledge, this is the first time that LFA has been applied to additive overlapping Vanka relaxation.

We consider the Stokes equations as a model problem, with both $P_2-P_1$  and $Q_2-Q_1$ discretizations.  We propose two constructions of the patches for Vanka relaxation, and two approaches to determine relaxation weights. It is shown that using weighting based on patch  geometry outperforms a simpler approach. We also find that using small patches with low-degree Chebyshev iterations leads to more efficient multigrid algorithms than with bigger patches or more relaxation steps per iteration, when cost per sweep is accounted for. Although there are no general rules to facilitate the choice of patches or weights, taking advantage of LFA, we can optimize the weights. For validation, our numerical tests are  implemented  using  Firedrake\cite{rathgeber2017firedrake} and PETSc.\cite{balay2018petsc,kirby2018solver} Numerical experiments are shown to match the LFA predictions for both periodic and Dirichlet boundary conditions. We observe that performance is less sensitive to overestimates of the weights for relaxation schemes considered here, which has also been seen in other works. \cite{JAdler_etal_2015b,adler2016monolithic} Last but not least, we compare the cost and performance of relaxation schemes considered here.

 This paper is organized as follows. In Section \ref{sec:FEM-Vanka}, we introduce the $P_2-P_1$ and $Q_2-Q_1$ discretizations considered here and the multigrid framework  with additive Vanka relaxation for the Stokes equations. In Section \ref{sec:LFA-intro}, we first give an introduction to LFA, then  propose an LFA for the Stokes equations with additive Vanka relaxation. In Section \ref{sec:Vanka-patch}, two types of overlapping patches are considered, and we validate the LFA predictions with two-grid experiments. Conclusions and remarks are given in Section \ref{sec:concl-FW}.

\section{Discretization and solution of the Stokes equations}\label{sec:FEM-Vanka}

\subsection{Mixed finite-element discretization of the Stokes
  equations}
In this paper, we consider the Stokes equations,
\begin{eqnarray}
  -\Delta\vec{u}+\nabla p&=&\vec{f},\label{Stokes-P2P1}\\
  -\nabla\cdot \vec{u}&=&0,\nonumber
\end{eqnarray}
where $\vec{u}$ is the velocity vector, $p$ is the  scalar pressure of a viscous fluid, and $\vec{f}$ represents a (known) forcing term, together with suitable boundary conditions.

The natural  finite-element approximation of Problem \eqref{Stokes-P2P1} when coupled with Dirichlet boundary conditions on $\vec{u}$ on some portion of the domain boundary is:
Find $\vec{u}_h\in \mathcal{X}^h$ and $p_h\in \mathcal{H}^{h}$  such that
\begin{equation}\label{mixed-FEM-form-P2P1}
a(\vec{u}_h,\vec{v}_h) + b(p_h,\vec{v}_h)+b(q_h,\vec{u}_h) =g(\vec{v}_h),\,\, {\rm for\,\, all}\, \vec{v}_h\in\mathcal{X}_{0}^h \,\,{\rm and}\,\, q_h\in\mathcal{H}^{h},
\end{equation}
where
\begin{eqnarray*}
a(\vec{u}_h,\vec{v}_h)& =& \int_{\Omega}\nabla \vec{u}_h:\nabla \vec{v}_h,\,\,\, b(p_h,\vec{v}_h)=-\int_{\Omega} p_h \nabla \cdot \vec{v}_h,\nonumber\\
 g(\vec{v}_h) &=&\int_{\Omega} \vec{f}_h\cdot \vec{v}_h,
\end{eqnarray*}
and $\mathcal{X}^h\subset H^{1}(\Omega)$, $\mathcal{H}^{h}\subset L_2(\Omega)$ are finite-element spaces. Here, $\mathcal{X}_{0}^h$ satisfies homogeneous Dirichlet boundary conditions in place of any inhomogenous essential boundary conditions on $\mathcal{X}^h$. Problem \eqref{mixed-FEM-form-P2P1} has a unique solution only when $\mathcal{X}^h$ and $\mathcal{H}^{h}$ satisfy an inf-sup condition. \cite{elman2006finite,MR2373954,brezzi1988stabilized,dohrmann2004stabilized}

\begin{remark}
If considering an outflow
   boundary condition, the stress-divergence form of the viscous
   term should be used instead\cite{limache2007}. The framework for LFA presented in this paper can easily be extended to this situation; however, the two-grid error-propagation operator (and its LFA representation) depends directly on the stencils of the discretized operators.  Thus, considering this form instead of \eqref{mixed-FEM-form-P2P1} may affect the optimal choice of parameters and the resulting performance of the two-grid method.
\end{remark}

\begin{figure}
\centering
\includegraphics[width=6.5cm,height=5.5cm]{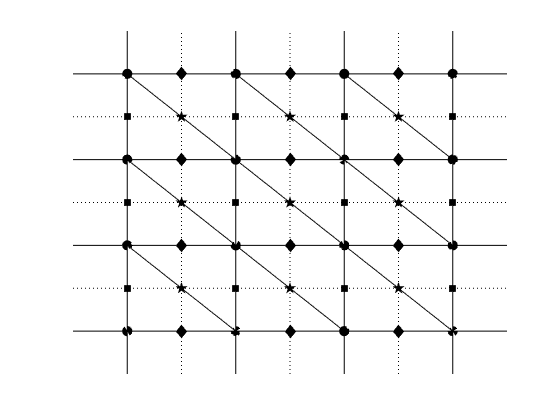}
\includegraphics[width=6.5cm,height=5.5cm]{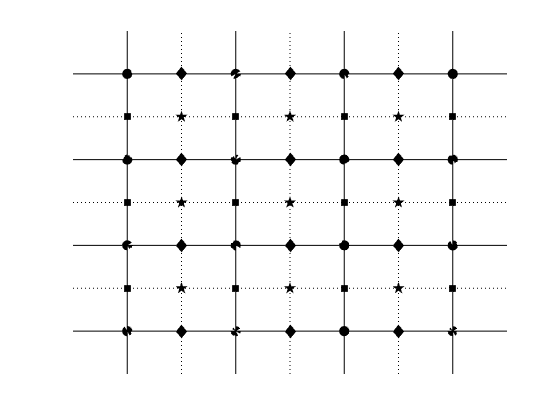}
\caption{Meshes and finite-element degrees of freedom (see definitions in \eqref{four-mesh-types}) , with \tikzcircle{3pt} denoting $N$-type  and $P_1/Q_1$ DoFs, $\blacklozenge$ denoting  $X$-type DoFs,  $\blacksquare$ denoting $Y$-type DoFs,  and $\bigstar$ denoting $C$-type DoFs. At left, $P_2-P_1$ discretization on triangles. At right, $Q_2-Q_1$ discretization on quadrilaterals.} \label{P2P1-Q2Q1-mesh}
\end{figure}

Here, we consider two types of stable finite-element methods for the Stokes equations.  First, we consider the stable mixed approximation  for structured meshes of triangular elements using continuous quadratic approximations
for the velocity components and continuous linear approximations for the pressure, the $P_2-P_1$ approximation.\cite{elman2006finite} Secondly, we consider the stable approximation for  rectangular meshes, using continuous biquadratic approximations
for the velocity components and continuous bilinear approximations for the pressure, the $Q_2-Q_1$ (Taylor--Hood) approximation. Both approximations  can be represented via nodal basis functions, as illustrated in Figure \ref{P2P1-Q2Q1-mesh}.

Discretizations of \eqref{Stokes-P2P1} typically lead to  linear systems of the form
\begin{equation}\label{saddle-structure-P2P1}
     Ky=\begin{pmatrix}
      A & B^{T}\\
     B & - C\\
    \end{pmatrix}
        \begin{pmatrix} {u }\\ {\rm p}\end{pmatrix}
  =\begin{pmatrix} {\rm f} \\ 0 \end{pmatrix}=b,
 \end{equation}
where $A$ corresponds to the discretized vector Laplacian, and $B$ is the negative of the discrete divergence operator. If the discretization is naturally unstable, then $C\neq 0$ is the stabilization matrix, otherwise $C=0$.  \cite{elman2006finite}  For the stable $P_2-P_1$ and $Q_2-Q_1$ finite-element discretizations considered here, we take $C=0$.

\subsection{Monolithic multigrid for the Stokes equations}

System \eqref{saddle-structure-P2P1} is of saddle--point type.  Here, we consider the numerical solution of \eqref{saddle-structure-P2P1} using a monolithic multigrid iteration applied to the full  system  collectively with a
suitable (coupled) relaxation method. As is typical for geometric multigrid, a relaxation technique is employed to quickly damp all oscillatory components of the error. Subsequently, a coarse-grid correction scheme, where a projected problem is solved on a coarser grid and the solution is
interpolated as an error correction to the fine-grid approximation, is used to damp the smooth components of the error.  In order to describe monolithic multigrid, assume we have two  meshes, with fine-grid meshsize $h$ and coarse-grid meshsize $H$ (often, $H=2h$, by doubling the meshsize in each spatial direction).

For a general nonsingular linear system, $K_h u_h=b_h$, we consider a stationary iteration as the relaxation scheme. Given an approximation, $M_h$, to $K_h$ that can be inverted easily, the approximate solution  is updated via the iteration
\begin{equation}\label{Stat-iteration}
  u^{j+1}_h = (I-M_h^{-1}K_h)u_h^{j} + M_h^{-1}b_h.
\end{equation}
The matrix $\mathcal{S}_h:=I-M_h^{-1}K_h$ is the error propagation operator for relaxation. With restriction and interpolation operators, $R_h$ and $P_h$, respectively, and coarse-grid matrix $K_H$,  the two-grid error propagation operator corresponding to the relaxation scheme in \eqref{Stat-iteration} can be written as
\begin{equation}\label{TG-Operator-P2P1}
  \boldsymbol{E}_h = \mathcal{S}_h^{\nu_2}\big(I-P_hK_H^{-1}R_hK_h\big)\mathcal{S}_h^{\nu_1},
\end{equation}
where  $I-P_hK_H^{-1}R_hK_h$ is called the coarse-grid correction operator.

The Jacobi, Gauss-Seidel, and Richardson schemes are often used for relaxation, particularly for discretizations of scalar PDEs.  For the restriction operator, $R_h$, there are many choices, which depend on the problem under consideration. Here, we focus on choices of $R_h$ tied to the mesh and the particular discretization scheme used to generate $K_h$. The coarse-grid operator, $K_H$ can be the Galerkin operator, $K_H= R_hK_hP_h$, or the natural rediscretization operator (or any other choice). The interpolation operator, $P_h$, is usually taken to be the conjugate transpose of $R_h$, with scaling depending on the discretization scheme and the dimension of the considered problem. For more details on the choice of multigrid components, see \cite{MR1807961,stuben1982multigrid,MR1156079}.

If we solve the coarse-grid problem recursively by the two-grid method, then we obtain a multigrid method. Over the past decades, a variety of types of multigrid methods have been developed, including $W,V,$ and $F$--cycles \cite{stuben1982multigrid}. In this paper, we focus on using additive Vanka-type relaxation in combination with a monolithic multigrid method to solve \eqref{saddle-structure-P2P1}. This means that $M_h$ is constructed by the Vanka approach, and updates both  components of  the solution to \eqref{saddle-structure-P2P1}  at the same time in the relaxation scheme. Only two-grid schemes are considered.

\subsubsection{Overlapping Schwarz relaxation}\label{Schwarz-alg}

Here, we present the multiplicative and additive Schwarz approaches to solve $K_hu_h=b_h.$  Let the degrees of freedom (DoFs) of $u_h$ be the set $\mathfrak{S}$, and $\mathfrak{S}_{i}, i=1,2,\cdots, N$, be  subsets of unknowns with $\mathfrak{S}=\bigcup_{i=1}^{N} \mathfrak{S}_i$. Let  $V_i$ be the restriction operator mapping  from vectors over the set  of all unknowns,  $\mathfrak{S}$, to  vectors whose unknowns consist of the DoFs in $\mathfrak{S}_{i}$. Then  $K_i=V_iK_hV_i^T$ is the restriction of $K_h$ to the $i$-th block of DoFs. Moreover, let $D_i={\rm diag}(d^{i}_1,d^{i}_2,\cdots,d^{i}_{m_i})$  for $i=1,\cdots,N$ be a diagonal weight matrix for each block $i$, where $m_i$ is the dimension of $K_i$. Then, the multiplicative and additive Schwarz iterations are presented in Algorithm \ref{MSM-alg}.

\begin{algorithm}
\caption{}\label{MSM-alg}
{\bf Multiplicative Schwarz iteration:}\\
 $u^{j,0}=u^{j-1,N}$. For $i=1,\cdots,N$,
  \begin{equation*}
      K_i \delta u_i=V_i(b_h-K_hu^{j,i-1}),
  \end{equation*}
  and
  \begin{equation*}
   u^{j,i}= u^{j,i-1} +V_i^{T}D_i\delta u_i.
  \end{equation*}

 {\bf Additive Schwarz iteration:}\\

For $i=1,\cdots,N$,
  \begin{equation*}
      K_i \delta u_i=V_i(b_h-K_hu^j),
  \end{equation*}
and
  \begin{equation*}
    u^{j+1} = u^j +\sum_{i=1}^{N}V_i^{T}D_i \delta u_i.
  \end{equation*}
\end{algorithm}

The error-propagation operator for the multiplicative Schwarz procedure  can be written as
\begin{equation*}
  S_m=\prod_{i=1}^{N}\big(I-V_i^{T}D_i K_i^{-1}V_iK_h\big)
\end{equation*}
 and, for the additive Schwarz procedure,  it is
\begin{equation*}
  S_a=I-M_h^{-1}K_h,
\end{equation*}
where \begin{equation}\label{ASM-precondition}
  M_h^{-1} = \sum_{i=1}^{N}V_i^{T}D_i K_i^{-1}V_i.
\end{equation}
More details about this algebraic viewpoint on the multiplicative and additive Schwarz iterations  can be found, for example, in the work of Saad. \cite{saad2003iterative}  Overlapping  multiplicative Schwarz approaches have been used as  the relaxation scheme for the Stokes equations,\cite{JAdler_etal_2015b,SPVanka_1986a,VJohn_LTobiska_2000a,VJohn_GMatthies_2001a,SManservisi_2006a,schoberl2003schwarz,rodrigo2016local,SPMacLachlan_CWOosterlee_2011a} which we refer to as  multiplicative Vanka relaxation. Extending this, we will refer  to such overlapping  additive Schwarz approaches as additive Vanka relaxation. In this paper, we focus on using additive Vanka as a relaxation scheme within monolithic multigrid for the  Stokes equations.
Some key questions in doing this are
\begin{enumerate}
\item How should the subsets $\mathfrak{S}_{i}$ be chosen?
\item How should $D_i$ be chosen?
\end{enumerate}
In what follows, we consider uniform meshes for the domain, $\Omega =[0,1]^2$. We will use the pressure DoFs to ``seed'' the sets, $\mathfrak{S}_i$, so that (away from the boundary) all sets will have the same structure and size $m$. In this paper, we use local Fourier analysis to guide the choice of $\mathfrak{S}_i$ and other aspects of the relaxation scheme.

\revise{LFA for multiplicative Vanka-type relaxation for the $Q_2-Q_1$ discretization of the Stokes equations was first performed by MacLachlan and Oosterlee \cite{SPMacLachlan_CWOosterlee_2011a}, building on early work by Molenaar\cite{molenaar1991two} and Sivaloganathan\cite{sivaloganathan1991use}.  This was extended to the $P_2-P_1$ discretization on triangular meshes by Rodrigo et al.\cite{rodrigo2016local}.  In the multiplicative context, the overlap requires careful treatment within LFA, leading to Fourier representations of so-called ``stages'' of error within the algorithm.  For the additive case considered here, the analysis is much simpler, akin to that of multicoloured relaxation\cite{boonen2008local, KahlKintscher_2018a}.  While we do not consider LFA for multiplicative Vanka here, we do offer numerical comparisons of performance between the additive and multiplicative variants below, when using the same subsets, $\mathfrak{S}_{i}$, to define the relaxation schemes.}

\revise{
  The use of additive Schwarz-type iterations as relaxation schemes in multigrid methods for saddle-point problems was previously analyzed by Sch\"oberl and Zulehner\cite{schoberl2003schwarz}. The authors show that, under suitable conditions, the iteration can be interpreted as a symmetric inexact Uzawa method. The advantage of this analysis compared to LFA is that it applies to unstructured grids, although at the cost of giving less quantitative insight. They also propose a general strategy for constructing patches (gathering the velocity degrees of freedom connected to given pressure degrees of freedom) and weights, and apply this to the Crouzeix--Raviart discretization of the Stokes equations. This patch strategy yields the so-called inclusive patches considered subsequently in Figure \ref{patch-P2P1-plot}, although with different weights (as these are tuned by LFA in this work).}

\revise{
  More generally, both additive and multiplicative Vanka relaxation can be viewed as domain decomposition methods.  Thus, the algorithms considered here could also be analyzed from the perspective of one- and two-level domain decomposition approaches, albeit with smaller than typical subdomains. The work of Szyld and Frommer\cite{MR1712686} seems most relevant to the relaxation schemes considered here, analysing the convergence of overlapping additive Schwarz, but in the case of $M-$matrices.  Similar analysis exists for overlapping multiplicative Schwarz, by Benzi et al.\cite{MR1865505}.  To our knowledge, such approaches have yet to be applied to Vanka-type relaxation.
}

\section{Local Fourier Analysis}\label{sec:LFA-intro}

\subsection{Definitions and notations}
 We first introduce some terminology of LFA.\cite{MR1807961,wienands2004practical}
We consider the following two-dimensional infinite uniform grids, $\mathbf{G}_h=\bigcup_{j=1}^{4}\mathbf{G}_h^j$, where
\begin{equation}\label{four-mesh-types}
  \mathbf{G}^{j}_{h}=\big\{\boldsymbol{x}^{j}:=(x^j_1,x^j_2)=(k_{1},k_{2})h+\delta^{j},(k_1,k_2)\in \mathbb{Z}^2\big\},
\end{equation}
with
\begin{equation*}
\delta^{j}=\left\{
  \begin{aligned}
    &(0,0) &\text{if}\quad j=1,\\
    &(h/2,0)  &\text{if} \quad j=2, \\
    &(0,h/2)  &\text{if} \quad j=3, \\
    &(h/2,h/2) &\text{if}\quad j=4.\\
  \end{aligned}
               \right.
\end{equation*}
We refer to $\boldsymbol{G}^{1}_{h}, \boldsymbol{G}^{2}_{h},\boldsymbol{G}^{3}_{h},$ and  $\boldsymbol{G}^{4}_{h}$ as the $N$-, $X$-, $Y$-, and $C$-type points on the grid $\boldsymbol{G}_h$, respectively, see Figure \ref{P2P1-Q2Q1-mesh}. The coarse grids, $\mathbf{G}^j_{2h}$, are defined similarly. Note that in much of the literature,  LFA is applied to  discretizations on $\mathbf{G}_h^1$. Here, we consider the more general case as needed for $P_2$ and $Q_2$ finite elements.

Let $L_h$ be a scalar Toeplitz operator defined by its stencil acting on $l^2(\mathbf{G}^{j}_{h})$ as follows,
\begin{eqnarray}\label{defi-symbol-P2P1}
  L_{h}  &\overset{\wedge}{=}& [s_{\boldsymbol{\kappa}}]_{h} \,\,(\boldsymbol{\kappa}=(\kappa_{1},\kappa_{2})\in \boldsymbol{V}); \,
  L_{h}w_{h}(\boldsymbol{x}^j)=\sum_{\boldsymbol{\kappa}\in\boldsymbol{V}}s_{\boldsymbol{\kappa}}w_{h}(\boldsymbol{x}^j+\boldsymbol{\kappa}h),
  \end{eqnarray}
with constant coefficients $s_{\boldsymbol{\kappa}}\in \mathbb{R} \,(\textrm{or} \,\,\mathbb{C})$, where $w_{h}(\boldsymbol{x}^j)$ is
a function in $l^2(\mathbf{G}^j_{h})$. Here, $\boldsymbol{V}\subset \mathbb{Z}^2$ is a finite index set. Because $L_h$ is formally diagonalized by the Fourier modes $\varphi(\boldsymbol{\theta},\boldsymbol{x}^j)= e^{\iota\boldsymbol{\theta}\cdot\boldsymbol{x}^j/\boldsymbol{h}}=e^{\iota\theta_1x^j_1/h}e^{\iota \theta_2x^j_2/h}$, where $\boldsymbol{\theta}=(\theta_1,\theta_2)$ and $\iota^2=-1$, we use $\varphi(\boldsymbol{\theta},\boldsymbol{x}^j)$ as a Fourier basis with $\boldsymbol{\theta}\in \big[-\frac{\pi}{2},\frac{3\pi}{2}\big)^{2}$ (or any pair of  intervals with length $2\pi$).

For smoothing and two-grid analysis, we have to distinguish high and low
frequency components on $\mathbf{G}^j_h$ with respect to $\mathbf{G}^j_{2h}$.\cite{MR1807961} Note that for any $\bm{\theta'} \in \Big[-\frac{\pi}{2},\frac{\pi}{2}\Big)^2$,
\begin{equation}\label{distiguish-LH}
 \varphi(\boldsymbol{\theta},\bm{x}^j) =\varphi(\boldsymbol{\theta'},\bm {x}^j)\,\, \text{for}\, \bm{x}^j \in \mathbf{G}^j_{2h},
\end{equation}
if and only if $\bm{\theta}=\bm{\theta'}({\rm mod}\,\, \pi)$. This means that only those frequency components, $\varphi(\boldsymbol{\theta},\cdot)$,  with $\bm {\theta}\in \Big[-\frac{\pi}{2},\frac{\pi}{2}\Big)^2$ are distinguishable on $\mathbf{G}^j_{2h}$. Thus, high and low frequencies for standard coarsening ($H=2h$) are given by
\begin{equation*}
  \boldsymbol{\theta}\in T^{{\rm low}} =\left[-\frac{\pi}{2},\frac{\pi}{2}\right)^{2}, \, \boldsymbol{\theta}\in T^{{\rm high}} =\displaystyle \left[-\frac{\pi}{2},\frac{3\pi}{2}\right)^{2} \bigg\backslash \left[-\frac{\pi}{2},\frac{\pi}{2}\right)^{2}.
\end{equation*}

\begin{definition}\label{formulation-symbol-P2P1}
  We call $\widetilde{L}_{h}(\boldsymbol{\theta})=\displaystyle\sum_{\boldsymbol{\kappa}\in\boldsymbol{V}}s_{\boldsymbol{\kappa}}e^{\iota \boldsymbol{\theta}\cdot\boldsymbol{\kappa}}$ the symbol of $L_{h}$.
\end{definition}
 Note that for all  functions $\varphi(\boldsymbol{\theta},\boldsymbol{x}^j)$,
\begin{equation*}
   L_{h}\varphi(\boldsymbol{\theta},\boldsymbol{x}^j)= \widetilde{L}_{h} (\boldsymbol{\theta})\varphi(\boldsymbol{\theta},\boldsymbol{x}^j).
\end{equation*}

For a relaxation scheme, represented by matrix $M_h$ for operator $L_h$, the error-propagation operator for relaxation can be written as
\begin{equation*}
\mathcal{S}_h(\bm p)=I- M_h^{-1}(\bm p)L_h,
\end{equation*}
where $\boldsymbol{p}$ represents parameters within $M_h$. A typical relaxation scheme often reduces high-frequency error components quickly, but is slow to reduce low-frequency errors. Thus, it is natural to define the smoothing factor as follows.

\begin{definition}\label{err-prop-LFA-P2P1}
 The error-propagation symbol, $\widetilde{\mathcal{S}}_{h}(\boldsymbol{\theta},\boldsymbol{p})$, for relaxation scheme $\mathcal{S}_{h}(\bm p)$ on the infinite grid  $\mathbf{G}_{h}$ satisfies
\begin{equation*}
  \mathcal{S}_{h}(\bm p)\varphi(\boldsymbol{\theta},\boldsymbol{x})=\widetilde{\mathcal{S}}_{h}(\bm \theta,\bm p)\varphi(\boldsymbol{\theta},\boldsymbol{x}), \,\,\boldsymbol{\theta}\in \bigg[-\frac{\pi}{2},\frac{3\pi}{2}\bigg)^{2},
\end{equation*}
for all $\varphi(\boldsymbol{\theta},\boldsymbol{x})$, and the corresponding smoothing factor for $\mathcal{S}_{h}(\bm p)$ is given by
\begin{equation*}
  \mu_{{\rm loc}}=\mu_{{\rm loc}}\big(\mathcal{S}_{h}(\bm p)\big)=\max_{\boldsymbol{\theta}\in T^{{
  \rm high}}}\Big\{\big|\widetilde{\mathcal{S}}_{h}(\boldsymbol{\theta},\bm p)\big| \,\,\Big\}.
\end{equation*}
\end{definition}

In many cases, the LFA smoothing factor offers a good prediction of actual two-grid performance, and we can optimize the smoothing factor with respect to the parameters, $\bm p$, to obtain an efficient algorithm.  However, this is generally not true for higher-order finite-element approximations.\cite{HMFEMStokes,HM2018LFALaplace,SPMacLachlan_CWOosterlee_2011a} Thus, we next introduce two-grid LFA, which still offers good predictions of performance in this setting.

\begin{remark}
In many applications of LFA, we consider a system operator rather than the discretization of a scalar PDE, and $\mathcal{S}_h$ is a block smoother. However, the definition of symbol and smoothing factor presented here can be extended to a system easily.  Details on these extensions will be presented in section \ref{LFA-stokes}.
\end{remark}

\subsection{Two-grid LFA}\label{subsec-two-LFA}
 In general,  LFA smoothing analysis  gives a good prediction for the actual multigrid performance, under the assumption that we have an ``ideal'' coarse-grid-correction operator that annihilates low-frequency error components and leaves  high-frequency components unchanged. However, in our setting, this assumption about ideal coarse-grid correction (CGC) does not hold (due to the discretization \cite{HM2018LFALaplace}), but  the two-grid LFA convergence factor still offers useful predictions.

To apply LFA to the two-grid operator, \eqref{TG-Operator-P2P1},  and calculate  the two-grid convergence factor, we need to analyse how the operators $K_h, P_h, R_h,$ and $\mathcal{S}_h$  act on the Fourier components $\varphi(\boldsymbol{\theta},\boldsymbol{x}^j)$.  From \eqref{distiguish-LH}, we know that values of $\varphi(\boldsymbol{\theta},\boldsymbol{x}^j)$ coincide on $\mathbf{G}^j_{2h}$ for four values of $\boldsymbol{\theta}$, known as harmonic frequencies. Let
\begin{eqnarray*}
\boldsymbol{\alpha}&=&(\alpha_1,\alpha_2)\in\big\{(0,0),(1,0),(0,1),(1,1)\big\},\\
\boldsymbol{\theta}^{\boldsymbol{\alpha}}&=&(\theta_1^{\alpha_1},\theta_2^{\alpha_2})=\boldsymbol{\theta}+\pi\cdot\boldsymbol{\alpha},\,\,
\boldsymbol{\theta}:=\boldsymbol{\theta}^{00}\in T^{{\rm low}}.
\end{eqnarray*}
For a given $\bm \theta\in T^{\rm low}$, we define the four-dimensional harmonic space
\begin{equation*}
  \mathcal{F}(\bm \theta)={\rm span}\Big\{ \varphi(\boldsymbol{\theta^{\alpha}},\cdot): \boldsymbol{\alpha}\in\big\{(0,0),(1,0),(0,1),(1,1)\big\}\Big\}.
\end{equation*}
Under standard assumptions, the space $\mathcal{F}(\bm \theta)$ is invariant under the two-grid operator $\boldsymbol{E}_h$. \cite{MR1807961,wienands2004practical} We use the ordering of $\boldsymbol{\alpha}=(0,0),(1,0),(0,1),(1,1)$ for the four harmonics in the following, although, as with any invariant subspace, the ordering of the basis elements is unimportant.

Inserting the representations of $\mathcal{S}_h, K_h, K_{H}, P_h, R_h$ into \eqref{TG-Operator-P2P1}, we obtain the Fourier representation of two-grid error-propagation operator as
\begin{equation*}
 \widetilde{ \boldsymbol{E}}_h(\boldsymbol{\theta},p)= \widetilde{\boldsymbol {S}}^{\nu_2}_h(\boldsymbol{\theta},p)\big(I-\widetilde{\boldsymbol {P}}_h(\boldsymbol{\theta})(\widetilde{K}_{H}(2\boldsymbol{ \theta}))^{-1}\widetilde{\boldsymbol{ R}}_h(\boldsymbol{\theta})\widetilde{\boldsymbol{ \mathcal{K}}}_{h}(\boldsymbol{\theta})\big)\widetilde{\boldsymbol{ S}}^{\nu_1}_h(\boldsymbol{\theta},p),
\end{equation*}
where
\begin{eqnarray*}
\widetilde{\boldsymbol{\mathcal{K}}}_h(\boldsymbol{\theta})&=&\text{diag}\left\{\widetilde{K}_h(\boldsymbol{\theta}^{00}), \widetilde{K}_h(\boldsymbol{\theta}^{10}),\widetilde{K}_h(\boldsymbol{\theta}^{01}),
\widetilde{K}_h(\boldsymbol{\theta}^{11})\right\},\\
\widetilde{\boldsymbol{S}}_h(\boldsymbol{\theta},p)&=&\text{diag}\left\{\widetilde{\mathcal{S}}_h(\boldsymbol{\theta}^{00},p),
\widetilde{\mathcal{S}}_h(\boldsymbol{\theta}^{10},p),\widetilde{\mathcal{S}}_h(\boldsymbol{\theta}^{01},p),
\widetilde{\mathcal{S}}_h(\boldsymbol{\theta}^{11},p)\right\},\\
\widetilde{\boldsymbol{R}}_h(\boldsymbol{\theta})&=&\left(\widetilde{R}_h(\boldsymbol{\theta}^{00}),\widetilde{R}_h(\boldsymbol{\theta}^{10}),
\widetilde{R}_h(\boldsymbol{\theta}^{01}),\widetilde{R}_h(\boldsymbol{\theta}^{11}) \right),\\
\widetilde{\boldsymbol{P}}_h(\boldsymbol{\theta})&=&\left(\widetilde{P}_h(\boldsymbol{\theta}^{00});\widetilde{P}_h(\boldsymbol{\theta}^{10});
\widetilde{P}_h(\boldsymbol{\theta}^{01});\widetilde{P}_h(\boldsymbol{\theta}^{11}) \right),
\end{eqnarray*}
in which ${\rm diag}\{T_1,T_2,T_3,T_4\}$ stands for the block diagonal matrix with diagonal blocks, $T_1, T_2, T_3$, and $T_4$.

\begin{definition}
The asymptotic two-grid convergence factor, $\rho_{{\rm asp}}$, is defined as
\begin{equation}\label{real-TGM}
  \rho_{{\rm asp}} = {\rm sup}\bigg\{\rho\big(\widetilde{\boldsymbol{E}}_h(\boldsymbol{\theta},p)\big): \boldsymbol{\theta}\in T^{{\rm low}}\bigg\}.
\end{equation}
\end{definition}
In practical use, we typically consider a discrete form of $\rho_{\rm asp}$, denoted by $\rho$, resulting from sampling $\rho_{\rm asp}$ over only a finite set of frequencies. In many cases, $\rho$ provides a sharp prediction of actual two-grid performance. It is well known, for example, that  LFA gives the exact two-grid convergence factor  for problems with periodic boundary conditions.\cite{stuben1982multigrid}   The calculation of $\rho$ is much cheaper, however, than direct calculation of $\rho(E_h)$ from \eqref{TG-Operator-P2P1}.  More importantly, since $\rho$ is a function of the parameters, $\bm p$, arising from  the relaxation scheme (or the coarse-grid correction),  we can optimise $\rho$ to achieve an optimally  efficient algorithm. In our setting, such  parameters appear in the diagonal scaling matrices $D_i$ mentioned in Section \ref{Schwarz-alg}. One of our goals in this paper is to use LFA to optimise the two-grid convergence factor of multigrid when using such relaxation schemes.

Next, we provide  details on LFA for additive Vanka relaxation for the Stokes equations.  Considering practical use, we focus on the stable $P_2-P_1$ discretization, as  is  easily generated using  general-purpose FEM tools on simplicial meshes, such as Firedrake\cite{rathgeber2017firedrake} and FEniCS.\cite{LoggMardalEtAl2012a,AlnaesBlechta2015a}  \revise{While LFA for this discretization has previously been considered by Rodrigo et al.\cite{rodrigo2016local}, we make use of a different mesh construction and Fourier basis than proposed there.  Thus, in the following sections, we provide full details of the LFA for both the discrete operator and finite-element interpolation operators in this setting.}  We will also show  LFA predictions  for additive Vanka relaxation for the $Q_2-Q_1$ discretization. \revise{LFA for the $Q_2-Q_1$ discretization of the Stokes equations have been previously presented in the work of He and MacLachlan \cite{HMFEMStokes}, where the authors consider different relaxation schemes, but the operator stencil and grid-transfer operators are the same as those used here. Since the key ingredient of the multigrid method of interest here is the use of the additive Vanka relaxation scheme, we present full details of the LFA representation of the relaxation schemes for both $P_2-P_1$ and $Q_2-Q_1$ discretizations.}

\subsection{LFA for $P_2-P_1$}\label{LFA-stokes}
In what follows, we consider the discretized  Stokes equations, which read
\begin{equation}\label{block-system-stencil-P2P1}
  K_h=\begin{pmatrix}
      -\Delta_{h} & 0 &  (\partial_{x})_{h}\\
      0 & -\Delta_{h} & (\partial_{y})_{h} \\
      -(\partial_{x})_{h}  & -(\partial_{y})_{h} & 0
    \end{pmatrix}
    :=\begin{pmatrix}
       K_h^{1,1}&         0 &              K_h^{1,3}  \\
       0&         K_h^{2,2} &              K_h^{2,3} \\
       K_h^{3,1}&         K_h^{3,2} &              0 \\
    \end{pmatrix}.
\end{equation}

For the $P_2-P_1$ discretization, the degrees of freedom for velocity are located on $\mathbf{G}_h=\bigcup_{j=1}^4\mathbf{G}_h^{j}$, containing four types of meshpoints as shown at the left of Figure \ref{P2P1-Q2Q1-mesh}. The Laplace operator in \eqref{block-system-stencil-P2P1} is defined by its weak form restricted to the finite-element basis functions. Here, we use the standard nodal basis and, consequently, can write this in stencil form, extending \eqref{defi-symbol-P2P1}, with $\boldsymbol{V}$ taken to be a finite index set of values, $\boldsymbol{V}=V_{N}\bigcup V_{X}\bigcup V_{Y}\bigcup V_{C}$ with $V_{N}\subset\mathbb{Z}^2$, $V_{X}\subset\big\{(z_x+\frac{1}{2},z_y)|(z_x,z_y)\in\mathbb{Z}^2\big\}$, $V_{Y}\subset\big\{(z_x,z_y+\frac{1}{2})|(z_x,z_y)\in\mathbb{Z}^2\big\}$, and $V_{C}\subset\big\{(z_x+\frac{1}{2},z_y+\frac{1}{2})|(z_x,z_y)\in\mathbb{Z}^2\big\}$. With this, the (scalar) discrete Laplace operator is naturally treated as  a block operator,  and the Fourier representation of each block can be calculated based on Definition \ref{formulation-symbol-P2P1}, with the Fourier bases adapted to account for the staggering of the mesh points. Thus, the symbols of $K_h^{1,1}$ and  $K_h^{2,2}$  are  $4\times 4$ matrices.  Similarly to the Laplace operator, both terms in the gradient, $(\partial_x)_h$ and $(\partial_y)_h$, can be treated as  ($4\times 1$)-block operators. Then, the symbols of $K^{1,3}_h$ and $K^{2,3}_h$ are $4\times 1$ matrices, calculated based on Definition \ref{formulation-symbol-P2P1} adapted for the mesh staggering. The symbols of $K^{3,1}_h$ and $K^{3,2}_h$ are the conjugate transposes of those of $K^{1,3}_h$ and $K^{2,3}_h$, respectively.  Accordingly, $\widetilde{K}_{h}$ is a $9\times9$ matrix for the $P_2-P_1$ discretization.

We denote the symbols of the finite-element discretizations of the Stokes equations as
 \begin{equation*}
   \widetilde{K}_h(\theta_1,\theta_2)=\begin{pmatrix}
      \widetilde{ A}(\theta_1,\theta_2) & 0 &  \widetilde{B}_x^{T}(\theta_1,
      \theta_2)\\
      0 & \widetilde{A}(\theta_1,\theta_2) & \widetilde{B}_y^{T}(\theta_1,\theta_2) \\
      \widetilde{B}_x(\theta_1,\theta_2)  & \widetilde{B}_y(\theta_1,\theta_2) & 0
    \end{pmatrix}.
    \end{equation*}

\begin{figure}
\centering
\includegraphics[width=4.5cm,height=5.cm]{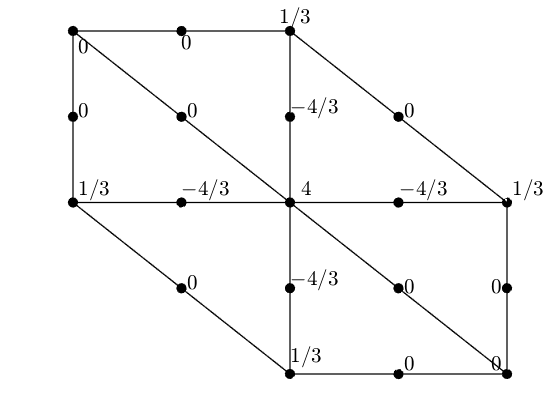}
\includegraphics[width=4.5cm,height=5.cm]{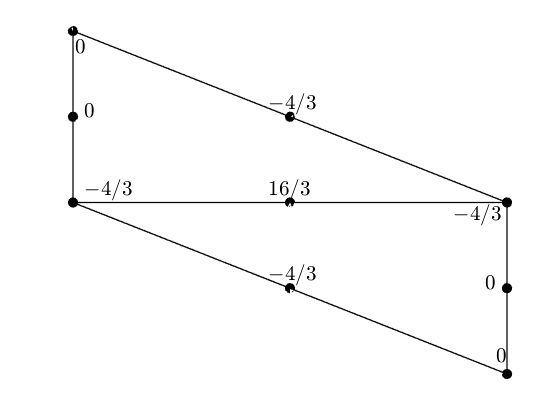}
\includegraphics[width=4.5cm,height=5.cm]{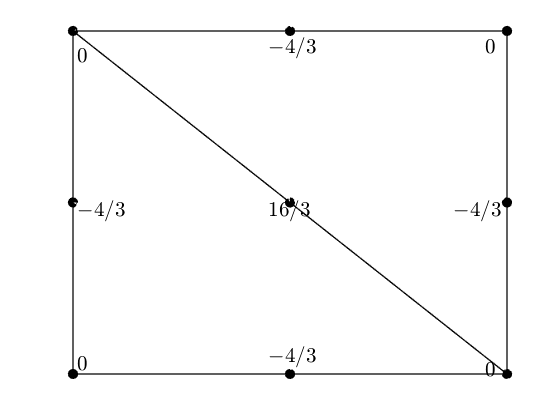}
\caption{Stencils for the $P_2$ finite-element discretization of the Laplace operator on a left triangular grid. Left: connections from DoF at a mesh node. Middle: connections from a horizontal edge. Right: connections from a diagonal mesh edge.} \label{N-X-stencil}
\end{figure}

Next, we discuss the stencils and  symbols for the operators in \eqref{block-system-stencil-P2P1}. For the Laplace operator, the stencil can be split into four types which  correspond to the $N$-, $X$-, $Y$-, and $C$-type points, shown  in Figure \ref{N-X-stencil}. For the $Y$-type, the stencil is a $90^{\circ}$ rotation of that of $X$-type, so we do not include it.

Rewriting the stencils shown in Figure \ref{N-X-stencil}, we can write the four  stencils of $K_h^{1,1}$  as follows,
\begin{eqnarray*}
 A_N&=&\frac{1}{3} \begin{bmatrix}
     &  & 1 &  &  \\
     &  & -4 &  &  \\
    1 & -4 & 12 & -4 &1 \\
     &  & -4 &  &  \\
     &  & 1 &  &
  \end{bmatrix},\,\, A_X= A_Y=A_C=\frac{1}{3} \begin{bmatrix}
         & -4 &   \\
       -4& 16 & -4  \\
         & -4 &
  \end{bmatrix}.
\end{eqnarray*}

Note that each stencil connects multiple types of meshpoints,  so we further  split each stencil into four  substencils based on the type of DoFs. Taking $A_N$ as an example, we see it connects three of the  four types of meshpoints. Thus,  $A_N$  can be written as  $A_N=\begin{pmatrix} A_{N,N} & A_{N,X}& A_{N,Y} &  A_{N,C} \end{pmatrix}$, where
 \begin{equation*}
   A_{N,N}=\frac{1}{3} \begin{bmatrix}
     &   1 &   \\
    1   & 12   &1 \\
     &   1   &
  \end{bmatrix},\,\,
  A_{N,X}=\frac{1}{3} \begin{bmatrix}
    -4   &    &-4
  \end{bmatrix},\,\,
  A_{N,Y}=\frac{1}{3} \begin{bmatrix}
        -4   \\
        -4
  \end{bmatrix},\,\,
  A_{N,C} =0.
 \end{equation*}

By standard calculation based on Definition \ref{formulation-symbol-P2P1}, we have
\begin{equation*}
\widetilde{A}_{N,N}= 4+ \frac{2}{3}\big(\cos\theta_1+\cos\theta_2\big),\,\, \widetilde{A}_{N,X} = -\frac{8}{3}\cos\frac{\theta_1}{2},\,\,
\widetilde{A}_{N,Y} = -\frac{8}{3}\cos\frac{\theta_2}{2},\,\, \widetilde{A}_{N,C} = 0.
\end{equation*}

Similarly, $\widetilde{A}_{X}, \widetilde{A}_{Y}, \widetilde{A}_{C}$ can be treated in this way. Thus, the symbol of $K_h^{1,1}$ and $K_h^{2,2}$  can be written as
\begin{equation*}
  \widetilde{A}(\theta_1,\theta_2)=\begin{pmatrix}
  \widetilde{A}_{N}\\
  \widetilde{A}_{X}\\
  \widetilde{A}_{Y}\\
  \widetilde{A}_{C}
  \end{pmatrix}
 = \begin{pmatrix}
  \widetilde{A}_{N,N} &  \widetilde{A}_{N,X}&  \widetilde{A}_{N,Y}&  \widetilde{A}_{N,C}\\
  \widetilde{A}_{X,N} &  \widetilde{A}_{X,X}&  \widetilde{A}_{X,Y}&  \widetilde{A}_{X,C}\\
  \widetilde{A}_{Y,N} &  \widetilde{A}_{Y,X}&  \widetilde{A}_{Y,Y}&  \widetilde{A}_{Y,C}\\
  \widetilde{A}_{C,N} &  \widetilde{A}_{C,X}&  \widetilde{A}_{C,Y}&  \widetilde{A}_{C,C}
  \end{pmatrix},
\end{equation*}
with $ \widetilde{A}(\theta_1,\theta_2)^T= \widetilde{A}(\theta_1,\theta_2)$.

 From Definition \ref{formulation-symbol-P2P1} and the stencils in Figure \ref{N-X-stencil}, we have  the following symbols,
\begin{eqnarray*}
\widetilde{A}_{X,X} &=& \frac{16}{3},\,\, \widetilde{A}_{X,Y} = 0,\,\,\widetilde{A}_{X,C} = -\frac{8}{3}\cos\frac{\theta_2}{2},\\
\widetilde{A}_{Y,Y} &=& \frac{16}{3},\,\, \widetilde{A}_{Y,C} = -\frac{8}{3}\cos\frac{\theta_1}{2},\,\,\widetilde{A}_{C,C} =\frac{16}{3}.
\end{eqnarray*}

\begin{figure}
\centering
\includegraphics[width=6.5cm,height=5.5cm]{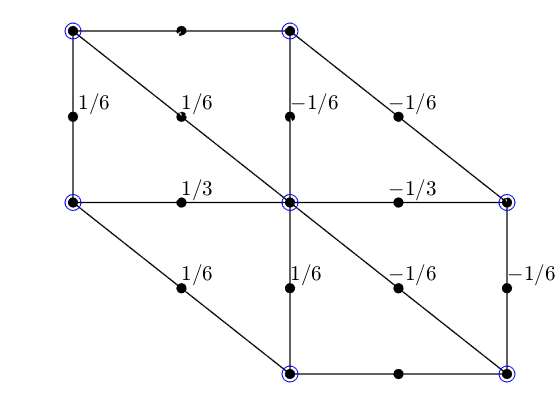}
\includegraphics[width=6.5cm,height=5.5cm]{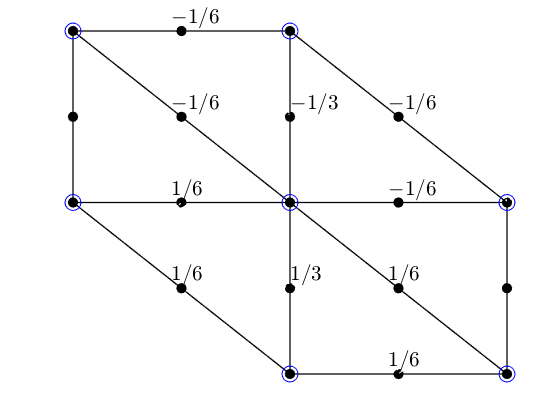}
\caption{Stencils for the $P_2-P_1$ finite-element discretization of the derivative operators on a left triangular grid. Left: $(\partial_x)_h$  stencil with a scaling $h$. Right: $(\partial_y)_h$  stencil with a scaling $h$. The (blue) circles at the center are the locations of pressure unknowns and the marked points without a given weight have value 0.} \label{BX-BY}
\end{figure}

Similarly to the stencil of the Laplacian operator, the stencils of $(\partial_x)_h$ and $(\partial_y)_h$  can be split into four types of substencil, respectively. Figure \ref{BX-BY} shows the stencils of the gradient, that is, the pressure-to-velocity unknowns ($N$-, $X$-, $Y$- and $C$-type) connections, for the pressure unknown located at the middle of the hexagon.

Thus, the stencil of $(\partial_{x})_h$ shown in Figure \ref{BX-BY} can be written as  $B_x^T= [B_{x,N};B_{x,X};B_{x,Y};B_{x,C}]$. However, here, we calculate the stencils of $-(\partial_{x})_h$ and its symbols given by

\begin{eqnarray*}
  B_{x,N}&=&0,\,\,
  \widetilde{B}_{x,N}(\theta_1,\theta_2)=0,\\
 B_{x,X}&=&\frac{h}{3}\begin{bmatrix}
   -1  & & 1
   \end{bmatrix},\,\,\,\,\quad
  \widetilde{B}_{x,X}(\theta_1,\theta_2)=\frac{2ih}{3}\sin\frac{\theta_1}{2},\\
  B_{x,Y}&=&\frac{h}{6}\begin{bmatrix}
   -1 & 1 & 0\\
   &  &\\
   0 & -1 &1\\
   \end{bmatrix},\,\,
  \widetilde{B}_{x,Y}(\theta_1,\theta_2)=\frac{ih}{3}\Bigg(\sin\frac{\theta_2}{2}+\sin\theta_1\cos\frac{\theta_2}{2}-\cos\theta_1\sin\frac{\theta_2}{2}\Bigg),\\
  B_{x,C}&=&\frac{h}{6}\begin{bmatrix}
   -1 & & 1\\
     &  &\\
   -1 & & 1
   \end{bmatrix},\,\,\,\,\quad
  \widetilde{B}_{x,C}(\theta_1,\theta_2)=\frac{2ih}{3}\sin\frac{\theta_1}{2}\cos\frac{\theta_2}{2},
\end{eqnarray*}
respectively.

Similarly to $\widetilde{B}_{x}(\theta_1,\theta_2)^T$, the symbol of the stencil of $(\partial_{y})_h$ can be written as
\begin{equation*}
\widetilde{B}_{y}(\theta_1,\theta_2)^{T} =
 [\widetilde {B}_{x,N}(\theta_2,\theta_1);\widetilde {B}_{x,Y}(\theta_2,\theta_1);\widetilde{ B}_{x,X}(\theta_2,\theta_1);
 \widetilde{ B}_{x,C}(\theta_2,\theta_1)].
 \end{equation*}

\begin{remark}
Rodrigo et al.~\cite{rodrigo2015local} presented a framework for LFA for edge-based
discretizations on triangular grids. They consider an expression for
the Fourier transform in a non-orthogonal coordinate system in space and frequency variables to adapt to arbitrary structured triangular  meshes.  This idea can be applied to the discretization considered here. However, for the triangle mesh considered here, it is not necessary to use non-orthogonal coordinates. In the framework of  LFA provided here, we use a different Fourier basis to be consistent with the different types of stencil located on different types of grid-points, which simplifies the calculation.

\end{remark}
\subsection{LFA representation of Grid-transfer operators}\label{RP-Operator}
Here, we use the standard finite-element interpolation operators and their transposes for restriction. In the following, we discuss the stencils and symbols of these restriction and interpolation operators.

 To derive symbols for the grid-transfer operators, we first consider an arbitrary restriction operator characterized by a constant coefficient
stencil $R_h\overset{\wedge}{=} [r_{\boldsymbol{\kappa}}]$. Then, an infinite grid function $w_{h}: \mathbf{G}_h^{1} \rightarrow \mathbb{R} \,(\textrm{or} \,\,\mathbb{C})$ is
transferred to the coarse grid, $\mathbf{ G}_{2h}^{1}$, in the following way
\begin{eqnarray*}
 ( R_h w_{h})(\boldsymbol{ x})&=&\sum_{\kappa\in{V}}r_{\boldsymbol{\kappa}}w_{h}(\boldsymbol{ x}+\boldsymbol{\kappa} h)\,\,(\boldsymbol{ x}\in \mathbf{G}_{2h}^{1}).
\end{eqnarray*}
Taking $w_h$ to be the Fourier mode, $\varphi(\boldsymbol{\theta^{\alpha}},\boldsymbol{x} )=e^{\iota\boldsymbol{\theta^{\alpha}}\cdot\boldsymbol{x}/\boldsymbol{h}}$, we have
\begin{eqnarray}\label{N-type-symbol}
 ( R_h \varphi(\boldsymbol{\theta^{\alpha}},\cdot ))(\boldsymbol{ x})&=& \widetilde{R}_h(\boldsymbol{\theta^\alpha}) \varphi_{2h}(2\boldsymbol{\theta}^{(0,0)},\boldsymbol{x})\,\,(\boldsymbol{ x}\in \mathbf{G}_{2h}^{1}),
 \end{eqnarray}
 with $\widetilde{R}_h(\boldsymbol{\theta^\alpha})=\displaystyle\sum_{\boldsymbol{\kappa} \in{V}}r_{\boldsymbol{\kappa}}e^{\iota\boldsymbol{\kappa}\cdot\boldsymbol{\theta^{\alpha}}}$, which is called the symbol of $R_h$. However, since we  consider discretizations on staggered meshes, where different ``types'' of variables  interact in the interpolation and
  restriction operators, the symbol definition for the restriction operator acting on $\mathbf{G}_{2h}^{j}$, where $j=2,3,4$, must be modified.

Similarly to the stencils of $K_h$,  the restriction operator for the components of velocity can also be decomposed based on the partitioning of the DoFs associated with the $N$-, $X$-,  $Y$-, and $C$-type meshpoints. Each of these restriction operators  connects between all four types of meshpoints, and we partition each restriction operator into four blocks based on the DoFs. For ${\bm x}\in \mathbf{G}_{2h}^{j}$,  where $j=1,\cdots,4$, by standard calculation,\cite{HM2018LFALaplace} \eqref{N-type-symbol} becomes
\begin{eqnarray*}
 ( R_h \varphi(\boldsymbol{\theta^{\alpha}},\cdot ))(\boldsymbol{ x})&=& \sum_{\kappa \in{V}}r_{\kappa}e^{\iota\boldsymbol{ \kappa}\cdot \boldsymbol{\theta^{\alpha}}}e^{\iota \boldsymbol{\alpha} \cdot \pi\boldsymbol{ x/h}} \varphi_{2h}(2\boldsymbol{\theta}^{(0,0)},\boldsymbol{x}).
 \end{eqnarray*}
 While $e^{\iota \boldsymbol{\alpha} \cdot \pi\boldsymbol{ x/h}}$ appears in the above formulation, it serves only to indicate which type of DoF $R_h$ is acting on, since
\begin{equation*}
e^{\iota \boldsymbol{\alpha} \cdot \pi\boldsymbol{ x/h}} = \left\{\begin{array}{cl}  1, & {\rm for}\,\,{\bm x} \in \mathbf{G}_{2h}^{1},\\
    (-1)^{\alpha_1}, & {\rm for}\,\,{\bm x}\in \mathbf{G}_{2h}^{2},\\
    (-1)^{\alpha_2}, & {\rm for}\,\,{\bm x} \in \mathbf{G}_{2h}^{3},\\
    (-1)^{\alpha_1}(-1)^{\alpha_2}, & {\rm for}\,\, {\bm x} \in \mathbf{G}_{2h}^{4}.
    \end{array}\right.
\end{equation*}
Thus, it is natural to give the following general  definition of a restriction symbol on a staggered mesh.
\begin{definition}\label{def-R-symbol-P2P1}
  We call $\widetilde{R}(\theta^\alpha)=\displaystyle\sum_{\kappa \in{V}}r_{\kappa}e^{\iota\boldsymbol{ \kappa}\cdot \boldsymbol{\theta^{\alpha}}}e^{\iota \boldsymbol{\alpha} \cdot \pi\boldsymbol{ x/h}}$ the restriction symbol of $R_h$.
\end{definition}
We emphasize that we must first split the restriction operator  into the different types of DoFs that it restricts from and to before we can apply  Definition \ref{def-R-symbol-P2P1}.

We first consider restriction to the $N$-type DoFs of a $P_2$ function,  which can be split into four blocks,
\begin{equation*}
R_{v,N}=\begin{bmatrix} R_{N,N}& R_{N,X} & R_{N,Y} & R_{N,C}\end{bmatrix}.
\end{equation*}
The $N$-to-$N$ connection is
\begin{equation*}
  R_{N,N}=\begin{bmatrix} 1\star\end{bmatrix},
\end{equation*}
where the $\star$ denotes the position (on the coarse grid) at which the discrete operator is applied. From Definition \ref{def-R-symbol-P2P1}, $\widetilde{R}_{N,N}=1$. The $X$-to-$N$ connections yield the stencil

\begin{equation*}
  R_{N,X}=\frac{1}{8}\begin{bmatrix}
  -1&     &     &     -1&           \\
   -1&   3&     \star&    3&        -1   \\
    &    -1&     &    &         -1
  \end{bmatrix}.
\end{equation*}
By standard calculation, we have
\begin{equation*}
  \widetilde{R}_{N,X}=\frac{1}{4}\Bigg(3\cos\frac{\theta_1}{2}-\cos\frac{3\theta_1}{2}-\cos\frac{\theta_1}{2}\cos\theta_2
  +\sin\frac{\theta_1}{2}\sin\theta_2-\cos\frac{3\theta_1}{2}\cos\theta_2-\sin\frac{3\theta_1}{2}\sin\theta_2\Bigg).
\end{equation*}
Similarly, the $Y$-to-$N$ connection has the stencil
\begin{equation*}
  R_{N,Y}=\frac{1}{8}\begin{bmatrix}
  -1&     -1&       \\
   &      3&    -1   \\
   &      \star & \\
   -1&    3&     \\
   &      -1&  -1
  \end{bmatrix},
\end{equation*}
with its symbol
\begin{equation*}
  \widetilde{R}_{N,Y}=\frac{1}{4}\Bigg(3\cos\frac{\theta_2}{2}-\cos\frac{3\theta_2}{2}-\cos\theta_1\cos\frac{\theta_2}{2}+\sin\theta_1\sin\frac{\theta_2}{2}
  -\cos\theta_1\cos\frac{3\theta_2}{2}-\sin\theta_1\sin\frac{3\theta_2}{2}\Bigg).
\end{equation*}
The $C$-to-$N$ connection has the stencil
\begin{equation*}
  R_{N,C}=\frac{1}{8}\begin{bmatrix}
   -1&      -1&   &   &      \\
   -1&       3&   &   &       \\
     &        & \star &   &\\
    &      &    &  3&   -1\\
    &      &    &  -1&   -1
  \end{bmatrix},
\end{equation*}

with its symbol
\begin{eqnarray*}
  \widetilde{R}_{N,C}&=&\frac{1}{8}\Bigg(3\cos\frac{\theta_1}{2}\cos\frac{\theta_2}{2}-3\sin\frac{\theta_1}{2}\sin\frac{\theta_2}{2} -
 \cos\frac{3\theta_1}{2}\cos\frac{3\theta_2}{2}-\sin\frac{3\theta_1}{2}\sin\frac{3\theta_2}{2}\\
 &&-\cos\frac{3\theta_1}{2}\cos\frac{\theta_2}{2}-\sin\frac{3\theta_1}{2}\sin\frac{\theta_2}{2}-
 \cos\frac{\theta_1}{2}\cos\frac{3\theta_2}{2}-\sin\frac{\theta_1}{2}\sin\frac{3\theta_2}{2}\Bigg).
\end{eqnarray*}

\begin{figure}
\centering
\includegraphics[width=6.5cm,height=5.5cm]{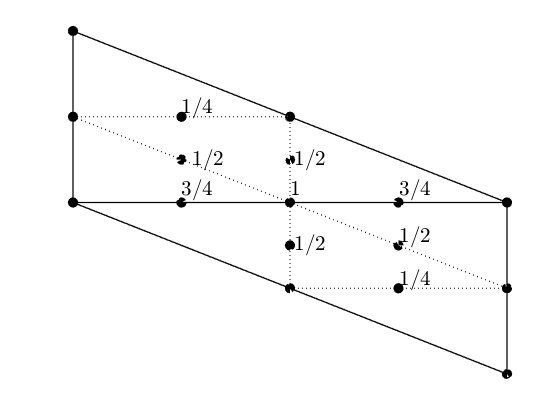}
\includegraphics[width=6.5cm,height=5.5cm]{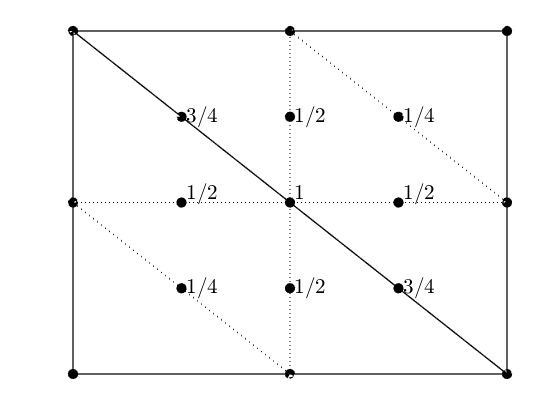}
\caption{Restriction stencils for the $P_2$ finite-element discretization of the Laplace operator on a left triangular grid. The solid triangles denote elements from the coarse-grid and the dashed triangles are those of  the fine-grid. Left:  $R_{v,X}$. Right: $R_{v,C}$.} \label{P2-restriction-plot}
\end{figure}

The weights for the $X$-type and $C$-type restrictions for $P_2$ are shown in Figure \ref{P2-restriction-plot}. For the $Y$-type DoFs, the restriction stencil is a $90^{\circ}$ rotation of that of $X$-type, so we do not include it.  We use the same decomposition for $R_{v,X},R_{v,Y}$ and $R_{v,C}$, and their symbols are listed in Table \ref{R-xyc-symbol}.

\begin{table}
\centering
\caption{Symbols of $R_{v,X},R_{v,Y}$ and $R_{v,C}$}
\begin{tabular}{l|l|l|l|l}
\hline
\hline
 *              &$R_{*,N}$      &$R_{*,X}$  & $ R_{*,Y}$    &$R_{*,C}$   \\
\hline
 X              &$1$              &$\frac{3}{2}\cos\theta_1+\frac{1}{2}(\cos\frac{\theta_1}{2}\cos\theta_2+\sin\frac{\theta_1}{2}\sin\theta_2)$    &$\cos\frac{\theta_2}{2}$                 & $\cos\frac{\theta_1}{2}\cos\frac{\theta_2}{2}+\sin\frac{\theta_1}{2}\sin\frac{\theta_2}{2}$                   \\
\hline
Y                     &$1$              &$\cos\frac{\theta_1}{2}$    &$\frac{3}{2}\cos\theta_2+\frac{1}{2}(\cos\theta_1\cos\frac{\theta_2}{2} +\sin\theta_1\sin\frac{\theta_2}{2})$                  & $\cos\frac{\theta_1}{2}\cos\frac{\theta_2}{2}+\sin\frac{\theta_1}{2}\sin\frac{\theta_2}{2}$                  \\
\hline
C                    &$1$              &$\cos\frac{\theta_1}{2}$   &$\cos\frac{\theta_2}{2}$                  &$2\cos\frac{\theta_1}{2}\cos\frac{\theta_2}{2}+\sin\frac{\theta_1}{2}\sin\frac{\theta_2}{2}$                   \\
\hline
\hline
\end{tabular}\label{R-xyc-symbol}
\end{table}

Following Definition \ref{def-R-symbol-P2P1} to account for staggering, we can write
 \begin{equation*}
         \widetilde{R}_{v}(\boldsymbol{\theta}^{00})=\begin{pmatrix}
         \widetilde{R}_{v,N}(\boldsymbol{\theta}^{00})\\
         \widetilde{R}_{v,X}(\boldsymbol{\theta}^{00})\\
         \widetilde{R}_{v,Y}(\boldsymbol{\theta}^{00})\\
         \widetilde{R}_{v,C}(\boldsymbol{\theta}^{00})
         \end{pmatrix},\,\,
          \widetilde{R}_{v}(\boldsymbol{\theta}^{10})=\begin{pmatrix}
         \widetilde{R}_{v,N}(\boldsymbol{\theta}^{10})\\
         -\widetilde{R}_{v,X}(\boldsymbol{\theta}^{10})\\
         \widetilde{R}_{v,Y}(\boldsymbol{\theta}^{10})\\
         -\widetilde{R}_{v,C}(\boldsymbol{\theta}^{10})
         \end{pmatrix},\,\,
          \widetilde{R}_{v}(\boldsymbol{\theta}^{01})=\begin{pmatrix}
         \widetilde{R}_{v,N}(\boldsymbol{\theta}^{01})\\
         \widetilde{R}_{v,X}(\boldsymbol{\theta}^{01})\\
         -\widetilde{R}_{v,Y}(\boldsymbol{\theta}^{01})\\
         -\widetilde{R}_{v,C}(\boldsymbol{\theta}^{01})
         \end{pmatrix},\,\,
          \widetilde{R}_{v}(\boldsymbol{\theta}^{11})=\begin{pmatrix}
         \widetilde{R}_{v,N}(\boldsymbol{\theta}^{11})\\
         -\widetilde{R}_{v,X}(\boldsymbol{\theta}^{11})\\
         -\widetilde{R}_{v,Y}(\boldsymbol{\theta}^{11})\\
         \widetilde{R}_{v,C}(\boldsymbol{\theta}^{11})
         \end{pmatrix}.
       \end{equation*}
Interpolation, $P_v$, is taken to be the transpose of $R_v$, with  symbol $\widetilde{P}_v(\bm {\theta^{\alpha}})=\frac{1}{4}\widetilde{R}^{T}_v(\bm{\theta^{\alpha}})$.

All of the pressure DoFs  are located on $\mathbf{G}_h^{1}$. Thus, the restriction operator for  pressure acts  only on  one type of grid point, the nodes of the mesh. The stencil of the restriction operator for pressure is given by
\begin{equation*}
  R_p=\frac{1}{2}\begin{bmatrix}
    1&   1        &   \\
    1&   2\star & 1   \\
     &   1        & 1
  \end{bmatrix},
\end{equation*}
with its symbol
\begin{equation*}
\widetilde{R}_p(\theta_1,\theta_2) = 1+ \cos\theta_1+\cos\theta_2+\cos\theta_1\cos\theta_2+\sin\theta_1\sin\theta_2.
\end{equation*}
 The interpolation, $P_p$, for pressure is taken to be the transpose of $R_p$, with symbol $\widetilde{ P}_p(\theta_1,\theta_2)=\frac{1}{4}\widetilde{R}(\theta_1,\theta_2)^{T}_p$.

Finally, the  restriction operator for the Stokes system can be written as
\begin{equation*}
R_h = \begin{bmatrix}
  R_v & 0 & 0 \\
  0 & R_v & 0 \\
  0 & 0 & R_p
\end{bmatrix},
\end{equation*}
with its symbol (over all four harmonics) being a $9\times 36$ matrix. For the interpolation operator, $P_h$, we consider  the transpose of $R_h$, and the symbols satisfy
\begin{equation*}
  \widetilde{\boldsymbol{ P}}_h(\bm \theta)=\frac{1}{4}\widetilde{\boldsymbol{ R}}_h^{T}(\bm \theta).
\end{equation*}
For the coarse-grid operator, $K_{H}$, we consider rediscretization, which is equivalent to Galerkin coarsening.

\subsection{LFA for additive Vanka relaxation}
When applying LFA to overlapping multiplicative Vanka relaxation, multiple Fourier representations are required for each variable in a given patch to account for the intermediate ``stages'' in the relaxation\cite{SPMacLachlan_CWOosterlee_2011a}. This is not needed in the additive case considered here.   Unlike classical relaxation, additive Vanka relaxation is a   block relaxation scheme. Since the variables on each block are updated at the same time, based on the same block system,  we can apply LFA ideas, but must modify the standard LFA to  include a block Fourier basis  to represent all the information in each block. By Fourier transformation the operator $M_h^{-1}$, defined on an infinite mesh, can be block diagonalized (by appropriately ordering the unknowns).

Recall $ M_h^{-1} = \sum_{i=1}^{N}V_i^{T}D_i K_i^{-1}V_i$. Under the Fourier ansatz, $D_i$, $V_i$  and $K_i$ have the same representation for all $i$. Note that since $D_i$ is a diagonal scaling  matrix, its Fourier representation is itself. Similarly, $V_i$ is a projector, and its symbol is itself.  Note that the representation of $K_i^{-1}$ is equal to the inverse of the representation of $K_i$. Thus, we only need to consider the representation of $K_i$.  Assume that the set of grid points  corresponding to the DoFs of $K_i$ is $\Xi_i=\big\{\boldsymbol{x}^{(i)}_1, \cdots, \boldsymbol{x}^{(i)}_m\big\}\subseteq \mathbf{G}_h$ and assume the ordering of the DoFs  in $K_i$ is consistent with the ordering of points in $\Xi_i$. Note that, due to the overlap between subdomains in the Schwarz relaxation, $\Xi_i$ might contain multiple points with the same Fourier representation in the symbol of \eqref{block-system-stencil-P2P1}, but we treat them separately in the representation of $K_i$, since these points  correspond to different unknowns in  $\mathfrak{S}_{i}$. Let
\begin{equation*}
  \varpi = {\rm span} \Big\{\psi_j(\boldsymbol{\theta}) = e^{\iota \boldsymbol{\theta}\cdot \boldsymbol{x_j^{(i)}}/h}\cdot \chi_j, j=1,\cdots,m\Big\},
\end{equation*}
where $\chi_j$ is an $m\times 1$ vector with only one nonzero element with value 1 located in the $j^{\mathrm{th}}$ position.

Let us consider the action of $K_i$ on this Fourier basis on $\Xi_i$, defining the symbol, $\widetilde{K}_{i}$, so that
\begin{equation}\label{Ki-acting}
  K_i\psi(\boldsymbol{\theta}) = \widetilde{K}_i \psi(\boldsymbol{\theta}),  \forall \psi \in \varpi.
\end{equation}
Let $\Phi$ be an $m\times m$ diagonal matrix with diagonal elements, $e^{\iota \boldsymbol{\theta}\cdot \boldsymbol{x}_j^{(i)}/h}$ for $j=1,2,\cdots,m$, and $\boldsymbol{x}^{(i)}_j \in \Xi_i$, and $\boldsymbol{\alpha} =\begin{pmatrix} \alpha_1,\alpha_2,\cdots,\alpha_m \end{pmatrix}^T$ be an arbitrary vector. Then, \eqref{Ki-acting} is equivalent to
\begin{equation}\label{FT-form}
  K_i \Sigma_{j=1}^m \alpha_j \psi_j  = K_i \Phi \boldsymbol{\alpha} = \Phi \widetilde{K}_i\boldsymbol{\alpha}.
\end{equation}
Thus, from \eqref{FT-form}, the Fourier representation of $K_i$ is
\begin{equation}\label{local-LFA}
  \widetilde{K}_i = \Phi^{-1}K_i\Phi = \Phi^TK_i\Phi.
\end{equation}
 Note that, under the Fourier ansatz, all of these matrices are independent of $i$ except for our construction of $\Phi$. However, due to the special structure of $\Phi$, we can use relative values of nodal position to replace the matrix $\Phi$  in \eqref{local-LFA} by a simple matrix (called the relative Fourier  matrix), scaled by $e^{\iota {\bm \theta}\cdot \boldsymbol{x}_s^{(i)}/h}$, for some $\boldsymbol{x}^{(i)}_s \in \Xi_i$. This is equivalent to considering the Fourier basis acting on a local offset, $\boldsymbol{x}^{(i)}_j:=\boldsymbol{x}^{(i)}_j-\boldsymbol{x}^{(i)}_s$ for some fixed point $\bm x_s^{(i)}$ in the patch, usually the pressure node.  As this scaling simplifies the calculation, we will use the relative Fourier  matrix in the rest of the paper.

\subsection{LFA for additive Vanka relaxation for the Stokes equations}

\begin{figure}
\centering
\includegraphics[width=6.5cm,height=5.5cm]{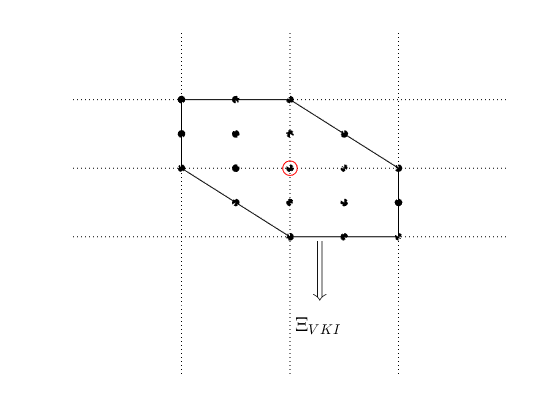}
\includegraphics[width=6.5cm,height=5.5cm]{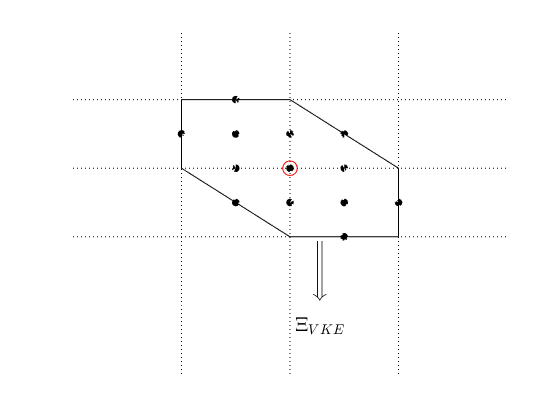}
\caption{Left: Vanka-inclusive patch for  the $P_2-P_1$  discretization. Right:  Vanka-exclusive patch for the $P_2-P_1$  discretization.} \label{patch-P2P1-plot}
\end{figure}

Now, we consider LFA for two types of Vanka relaxation schemes that differ only in the choice of the relaxation blocks. First, for each pressure DoF, we consider patches containing all velocity DoFs included in  a hexagon centred at the node associated with the pressure, see left of  Figure \ref{patch-P2P1-plot} for the $P_2-P_1$ discretization. This patch is the smallest one that contains all velocity DoFs with connections to this pressure DoF in the symbolic nonzero pattern of the matrix. For this reason, we refer to this patch construction as Vanka-inclusive (VKI).  The number of DoFs in this patch for $P_2-P_1$ is 39. To be specific, there are 7 $N$-type points, 4 points for $X$-,$Y$- and $C$-type, respectively, for each of the two components of the velocity, and 1 $N$-type point for the central pressure, see Table \ref{Patch-points-number}. Since the patches and submatrices  $K_i$ are the same, it is natural to consider $D_i$ to be the same for all $i$.  For $D_i$ in \eqref{ASM-precondition}, a simple  idea is to take $D_i$ to be the identity, which we refer to  as using {\it no weights}.  Another choice for $D_i$ is to take $d_{m_i}$ in $D_i$ to be the reciprocal of the number of patches that each type of DoF appears in. We refer to this as using {\it geometric weights} (or {\it natural weights}), and will denote this by VKIW in the results to follow.  To be specific, the weights for the velocity DoFs are $1/7,1/4,1/4,$ and $1/4$ for the $N$-, $X$-, $Y$- and $C$- type points, respectively. For pressure, the natural weight is 1.  Thus, $D_i$ is a $39\times 39$ matrix given by
\begin{equation*}
  D_i=  \begin{pmatrix}D_v & 0 &0\\ 0 & D_v  & 0\\ 0 & 0 &D_p \end{pmatrix}, \,\, {\rm with}\,\,D_v =\begin{pmatrix} D_N & 0 & 0 & 0\\0 & D_X & 0& 0 \\ 0 & 0& D_Y& 0\\ 0 & 0& 0& D_C\end{pmatrix},\,\, D_p=1,
\end{equation*}
where $D_N=\frac{1}{7}I_{7\times 7}$ and $D_X=D_Y=D_C=\frac{1}{4}I_{4\times 4}$. Since $\widetilde{K}_h$ is a $9 \times 9$ matrix,  the symbol of $V_i$ is a $39\times 9$ projection  matrix. Finally, $\widetilde{K}_i$ is a $39 \times 39$ matrix, as is $D_i$.

\begin{table}
\centering
\caption{The number (\#) of different type of DoFs for Vanka-inclusive and Vanka-exclusive patches for the $P_2-P_1$ and $Q_2-Q_1$ discretizations.}
\begin{tabular}{l||c|c||c|c}

        &\multicolumn{2} {c||} {$P_2-P_1$}  & \multicolumn{2} {c} {$Q_2-Q_1$} \\
\hline
 \backslashbox{\#}{patch}            &$\Xi_{VKI}$   & $\Xi_{VKE}$  &$\Xi_{VKI}$   &$\Xi_{VKE}$   \\
\hline
    \multicolumn{5} {c}  {each component of velocity} \\
\hline
N                 & 7      &1        &9       &1   \\
\hline
X                 &4      &4        &6       &6   \\
\hline
Y                 &4      &4        &6       &6   \\
\hline
C                 &4      &4        &4       &4   \\
\hline
    \multicolumn{5} {c}  {pressure} \\
\hline
N                 &1     &1        &1       &1   \\
\hline
\hline
Total            &39      &27      &51       &35  \\
\hline
\end{tabular}\label{Patch-points-number}
\end{table}

Now, we give more details about how to calculate the relative Fourier matrix, $\Phi$. For $P_2-P_1$, with patches $\Xi_{VKI}$ as shown  at the left of Figure \ref{patch-P2P1-plot}, we take $x^{(i)}_s$ be the node located directly below the central pressure DoF.  We use lexicographical order  to order the remaining  points in $\Xi_{VKI}$ for each type of point ($N$-,$X$-,$Y$-, and $C$- ordering), that is, from left to right and bottom to top. Then,
\begin{equation}\label{VKI-LFA}
\Phi(\theta_1,\theta_2) = \begin{pmatrix}\Phi_v & 0 &0\\ 0 & \Phi_v  & 0\\ 0 & 0 &\Phi_p \end{pmatrix},
\end{equation}
where
\begin{eqnarray*}
\Phi_n &=&[1,e^{\iota\theta_1}, e^{-\iota \theta_1+\iota\theta_2}, e^{\iota \theta_2}, e^{\iota \theta_1+\iota\theta_2}, e^{-\iota\theta_1+\iota2\theta_2},e^{\iota 2\theta_2}], \\
\Phi_x &=&[e^{\iota/2\theta_1},e^{-\iota/2\theta_1+\iota \theta_2}, e^{\iota/2\theta_1+\iota \theta_2},e^{-\iota/2\theta_1+\iota 2\theta_2}],\\
\Phi_y &=&[e^{\iota/2\theta_2}, e^{\iota\theta_1+\iota/2\theta_2},e^{-\iota\theta_1+\iota3/2\theta_2},e ^{\iota 3/2\theta_2}],\\
\Phi_c &=&[e^{-\iota/2\theta_1+\iota/2\theta_2},e^{\iota/2\theta_1+\iota/2\theta_2}, e^{-\iota/2\theta_1+\iota3/2\theta_2}, e^{\iota/2\theta_1+\iota3/2\theta_2}],\\
\Phi_{0} &=& [\Phi_n,\Phi_x,\Phi_y,\Phi_c],\\
\Phi_v & =& {\rm diag}(\Phi_{0}),\\
\Phi_p &=& e^{\iota \theta_2}.
\end{eqnarray*}

\begin{remark}
Note that matrix $\Phi$ in \eqref{VKI-LFA} is only a function of frequency, $\boldsymbol{\theta}=(\theta_1,\theta_2)$, and is independent of the meshsize, $h$, and the index $i$.
\end{remark}

For  $V_i$, mapping the global vector to the vector  on the patch,  we only need to account for the duplication that arises from the representation of the global Fourier basis of dimension 9 to the local block representation. Here, the order follows $N$-type, $X$-type, $Y$-type, then $C$-type, and velocity is first, then followed by the pressure.  The structure is \begin{equation*}
  \widetilde{V_i} = \begin{pmatrix} V_v & 0 & 0\\ 0 & V_v & 0 \\ 0 & 0 & V_p\end{pmatrix},
\end{equation*}
where  $V_p=1$  and $V_v$ is a $19\times 4$ matrix defined as follows
\begin{equation}\label{Vi-VKI}
  V_v = \begin{pmatrix} I_N & 0 & 0 & 0\\0 & I_X & 0& 0 \\ 0 & 0& I_Y& 0\\ 0 & 0& 0& I_C,
  \end{pmatrix}
\end{equation}
where $I_N$ is the $7\times 1$ vector of all ones, and $I_X=I_Y=I_C$ are the $4\times 1$ vector of all ones.

Considering the stencils in Figure \ref{BX-BY}, we note that of the adjacent nodal velocity DoFs, only the central one appears in the gradient operator.  This motivates construction of a second patch that  contains all of the DoFs in VKI except for the $N$-type velocity DoFs at the edges of the VKI patch, see the right of Figure \ref{patch-P2P1-plot}. We refer to this as Vanka-exclusive (VKE) when used with no weights ($D_i=I$). Table \ref{Patch-points-number} presents the number of the four types of DoFs for VKE, which contains $39-2\cdot6=27$ DoFs in total.  Another choice for $D_i$ is to use the natural weights, which  we refer to as Vanka-exclusive with natural weights (VKEW). Here, the weights are $1,1/4,1/4,1/4$ for $N$-, $X$-, $Y$- and $C$- type velocity DoFs, respectively, and 1 for pressure.  Similarly to VKI, we can calculate the matrix $\Phi$ whose size is $27\times 27$ following \eqref{VKI-LFA}, replacing $\Phi_n$ by  $\Phi_{n} =e^{\iota \theta_2}$.  For $\widetilde{V_i}$, there is only 1 $N$-type point for velocity, so we only need to modify $I_N$ to be the scalar 1 in \eqref{Vi-VKI}, then we obtain $\widetilde{V_i}$, which is a $27\times 9$ matrix. In each of these cases, the overall symbol for relaxation is a $9\times 9$ matrix.  Since  we consider four harmonics, following subsection \ref{subsec-two-LFA}, $\widetilde{ \boldsymbol{E}}_h$ is a $36\times 36$ matrix.  In the rest of this paper, we omit the subscript $h$ unless it is necessary to avoid confusion.

\begin{figure}
\centering
\includegraphics[width=6.5cm,height=5.5cm]{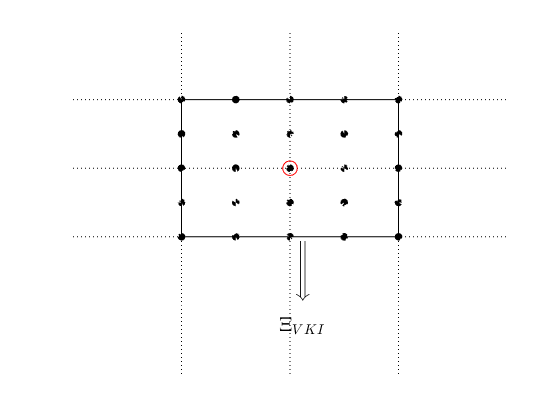}
\includegraphics[width=6.5cm,height=5.5cm]{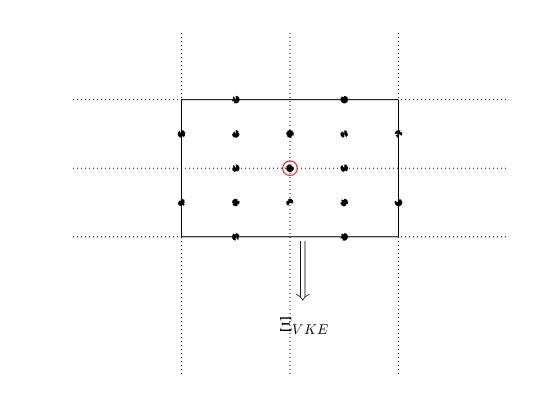}
\caption{Left: Vanka-inclusive patch for the $Q_2-Q_1$  discretization. Right: Vanka-exclusive patch for the $Q_2-Q_1$  discretization. } \label{patch-Q2Q1-plot}
\end{figure}

For the $Q_2-Q_1$ discretization of the Stokes equations, the patches of Vanka-inclusive and Vanka-exclusive are a little different than the structure with the $P_2-P_1$ discretizations. Figure \ref{patch-Q2Q1-plot} shows the Vanka-inclusive and Vanka-exclusive  patches for the $Q_2-Q_1$ discretization, and Table \ref{Patch-points-number} lists the details.   For VKI, the numbers of unknowns for the velocity are 9, 6, 6, 4 for $N$-, $X$-, $Y$-, and $C$-types, respectively.  In total, there are $2(9+6+6+4)+1=51$  DoFs in one patch.  To construct $V_{v}$, we only need to change $I_N$ in \eqref{Vi-VKI} to be the $9\times 1$ vector of all ones, and $I_X=I_Y$ to be the $6 \times 1$ vectors of all ones. As in the $P_2-P_1$ case, the $Q_2-Q_1$ gradient operator on a uniform mesh uses only the central nodal velocity DoF.  Thus, the second patch contains all of the DoFs in VKI except for the $N$-type points at the edges of the block, giving $51-2\cdot8=35$ DoFs. To obtain the representation of $V_i$ using VKE for $Q_2-Q_1$ from the construction of $V_v$ for VKI with the $Q_2-Q_1$ discretization, we only need to change $I_N$ to a scalar $1$ .

We again use the relative Fourier  matrix, $\Phi$, to transform, the block matrix, $K_i$.  We set $\boldsymbol{x}^{(i)}_s$ to be the lower-left corner at the block, giving,
\begin{equation}\label{VKI-LFA-Q2Q1}
\Phi(\theta_1,\theta_2) = \begin{pmatrix}\Phi_{Qv} & 0 &0\\ 0 & \Phi_{Qv}  & 0\\ 0 & 0 &\Phi_{Qp} \end{pmatrix},
\end{equation}
where
\begin{eqnarray*}
\Phi_{Qn} &=&[1,e^{\iota\theta_1}, e^{\iota2\theta_1}, e^{\iota\theta_2}, e^{\iota \theta_1+\iota \theta_2}, e^{\iota 2\theta_1+\iota\theta_2}, e^{\iota2\theta_2},e^{\iota\theta_1+\iota 2\theta_2},e^{\iota 2\theta_1+\iota 2\theta_2}], \\
\Phi_{Qx} &=&[e^{\iota/2\theta_1},e^{\iota3/2\theta_1},e^{\iota/2\theta_1+\iota \theta_2}, e^{\iota3/2\theta_1+\iota \theta_2},e^{\iota/2\theta_1+\iota 2\theta_2},e^{\iota3/2\theta_1+\iota 2\theta_2}],\\
\Phi_{Qy} &=&[e^{\iota/2\theta_2},e^{\iota\theta_1+\iota/2\theta_2}, e^{\iota2\theta_1+\iota/2\theta_2},e^{\iota3/2\theta_2},e ^{\iota\theta_1+\iota 3/2\theta_2},e ^{\iota2\theta_1+\iota 3/2\theta_2}],\\
\Phi_{Qc} &=&[e^{\iota/2\theta_1+\iota/2\theta_2},e^{\iota3/2\theta_1+\iota/2\theta_2}, e^{\iota/2\theta_1+\iota3/2\theta_2}, e^{\iota3/2\theta_1+\iota3/2\theta_2}],\\
\Phi_{Q0} &=& [\Phi_{Qn},\Phi_{Qx},\Phi_{Qy},\Phi_{Qc}],\\
\Phi_{Qv} & =& {\rm diag}(\Phi_{Q0}),\\
\Phi_{Qp} &=& e^{\iota\theta_1+\iota \theta_2}.
\end{eqnarray*}

For VKE with $Q_2-Q_1$, we truncate $\Phi_{Qn}$ to the scalar $ e^{\iota \theta_1+\iota \theta_2}$.

\section{Parameter choice and validation}\label{sec:Vanka-patch}

As a relaxation scheme, we  consider the  Chebyshev iteration \cite{saad2003iterative} on $K_h$ preconditioned with overlapping additive Vanka, with the two patches discussed before. We note that this is the natural extension of weighted Jacobi to block relaxation when considering more than a single relaxation sweep per level in a two-grid cycle. The key point to tuning the Chebyshev iteration is the choice of the lower and upper bounds for the interval that determines the Chebyshev polynomials. LFA is useful here, and we employ it to find optimal bounds for different degrees of Chebyshev polynomials. Relaxation using both no weights and natural weights will be considered.  As a special case, we also consider a simple  preconditioned Richardson relaxation, for the cases $\nu_1+\nu_2=1$ and $\nu_1+\nu_2=2$, where we again use LFA to help find the optimal weights. The goal of these experiments  is to use LFA to  determine ``best practices'' in terms of how
to choose patches and weights, with relaxation parameters optimized for these choices. Finally, we compare the cost and performance among these approaches.

In practice, the LFA two-grid convergence factors often exactly match the true convergence
factor of multigrid applied to a problem with periodic boundary conditions. \cite{MR1807961,stevenson1990validity} For the case of Dirichlet boundary conditions, a gap between the LFA predictions and the actual performance is sometimes observed. \cite{he2018local,HMFEMStokes,niestegge1990analysis} In order to see the influence of  boundary conditions
on multigrid performance, we present  data for both Dirichlet and periodic boundary conditions. Our tests are implemented using Firedrake and PETSc \cite{kirby2018solver} for both the $P_2-P_1$ and $Q_2-Q_1$ discretizations. The subproblems associated with $K_i$
  are solved by $LU$ decomposition. We focus on optimizing the interval used for Chebyshev relaxation with symmetric pre- and post-relaxation.
\subsection{Numerical results for the $P_2-P_1$ discretization}
\subsubsection{No weights}
First, we consider $D_i=I$ and the two different patches discussed above. The two-grid error iteration matrix is
\begin{equation}\label{Cheby-form-TG}
  E_{k} =  p_k(T)\mathcal{M}^{{\rm CGC}}p_k(T),
\end{equation}
where $T =M_h^{-1}K_h$ and $p_k$ is the Chebyshev polynomial with degree $k$ on a given interval.  Below, we take $\rho$  to be the  LFA prediction sampled at 32 equispaced points in each dimension of the Fourier domain. (When sampling the frequency space at more points, the predictions go up only slightly.) We take $\hat{\rho}$ to be the measured  convergence factor per iteration defined by $\hat{\rho}_j=\frac{||e_{j}||}{||e_{j-1}||}$, where $e_j=b_h-K_hu_j$, measured for problem \eqref{saddle-structure-P2P1} with zero right-hand side and random initial guess at the iteration when the residual norm is first below $10^{-150}$. Since the discretized  problem is singular, we project the approximate solution after each iteration to ensure that it remains orthogonal to the null space.

In order to find the optimal intervals for the Chebyshev iteration for different $k$, we use a brute force search with stepsize $0.1$ to optimize the intervals with LFA. The optimal interval, denoted as $[\alpha,\beta]$, and the corresponding convergence factors are given in the following tables. Tables \ref{No-TG-VKI} and \ref{No-TG-VKE} give the measured convergence   for periodic and Dirichlet boundary conditions versus  the LFA predictions for different $k$ and the corresponding optimal intervals with $h=\frac{1}{20},\frac{1}{40}$, and $\frac{1}{80}$. We see that the measured convergence factors  are largely  independent of the meshsize, $h$. A good match is generally seen between the LFA predictions and the measured convergence factors, with only  a tiny gap between the LFA prediction and the measured convergence factor. It is reasonable that such a gap exists for the Dirichlet case, since the LFA prediction is of the expected asymptotic convergence factor for the  problem with periodic boundary conditions. This might also suggest that extra work is needed to reduce this gap.\cite{niestegge1990analysis} Note that only for $k=2$ and VKI is this gap significant. We have also tested some other intervals for this case, seeing that the gap  between the LFA prediction and the measured results with Dirichlet boundary conditions is quite variable, with some parameter values coming close to the LFA prediction, but never achieving a good agreement.   For VKI with $k=3$ and periodic boundary conditions, the per iteration measured convergence factor oscillates, but the long-term  average still approaches that given by LFA.  For these cases, we present averaged convergence over the final seven iterations. These tables  indicate that when $k=1$, VKE is more effective than VKI; however, for $k>1$, the opposite occurs.

\begin{table}
\centering
\caption{Two-grid LFA predictions for Chebyshev Vanka relaxation and two-grid performance with periodic ($\hat{\rho}$) and Dirichlet ($\hat{\rho}^{\rm D}$) boundary conditions for the $P_2-P_1$ discretization for VKI. Optimal intervals for different $k$ with no weights. An alternating convergence pattern was observed, so the reported convergence factor is averaged over the final seven iterations before
convergence (denoted *). }
\begin{tabular}{l|l|l||l|l||l|l||l|l}
\hline
\multicolumn{3} {c||} {LFA predictions}  &\multicolumn{2} {c||} {$h=\frac{1}{20}$}   &\multicolumn{2} {c||} {$h=\frac{1}{40}$}  & \multicolumn{2} {c} {$h=\frac{1}{80}$} \\
\hline
   $k$  &$[\alpha,\beta]_{\rm VKI}$  &$\rho_{\rm VKI}$  & $\hat{\rho}_{\rm VKI}$  & $\hat{\rho}^{\rm D}_{\rm VKI}$   &$\hat{\rho}_{\rm VKI}$  & $\hat{\rho}^{\rm D}_{\rm VKI}$   &$\hat{\rho}_{\rm VKI}$  & $\hat{\rho}^{\rm D}_{\rm VKI}$ \\
\hline
 1                     &$[0.1,8.3]$              &0.672    &0.672   &0.699                &0.671  &0.699       &0.671  & 0.699    \\
\hline
 2                     &$[0.9,7.8]$              &0.295    &0.296   &0.585                &0.296  &0.585       &0.296   &0.585     \\
\hline
 3                     &$[0.9,7.9]$              &0.120    &0.107*   &0.148                &0.133*  &0.148       &0.143*    &0.148     \\
\hline
 4                     &$[1.4,7.2]$              &0.102    &0.102   &0.133                &0.102  &0.133       &0.102    &0.133     \\
\hline
 5                     &$[1.3,7.4]$              & 0.051   &0.050   &0.074                &0.049  &0.074       &0.047    &0.074   \\
\hline
\hline
\end{tabular}\label{No-TG-VKI}
\end{table}

\begin{figure}
\centering
\includegraphics[width=7.cm,height=5.5cm]{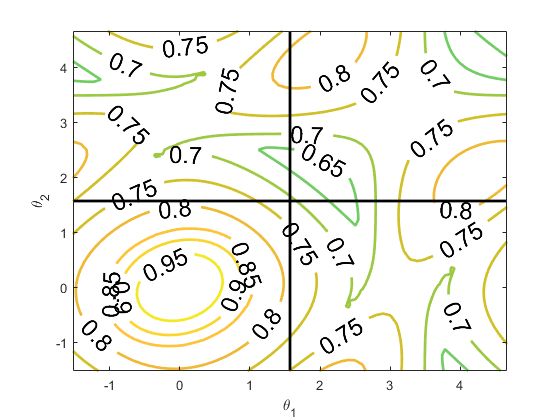}
\includegraphics[width=7.cm,height=5.5cm]{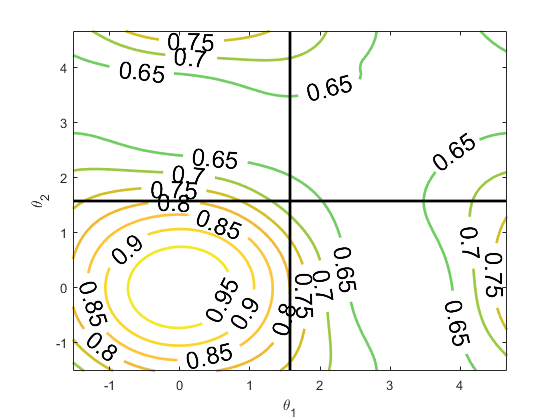}
\caption{The LFA  amplification factors for Vanka-relaxation with $k=1$ for the $P_2-P_1$ discretization. Left: VKI with $k=1$. Right: VKE with $k=1$.} \label{P2P1VK-SMOOTHING}
\end{figure}

\begin{figure}
\centering
\includegraphics[width=7.5cm,height=5.5cm]{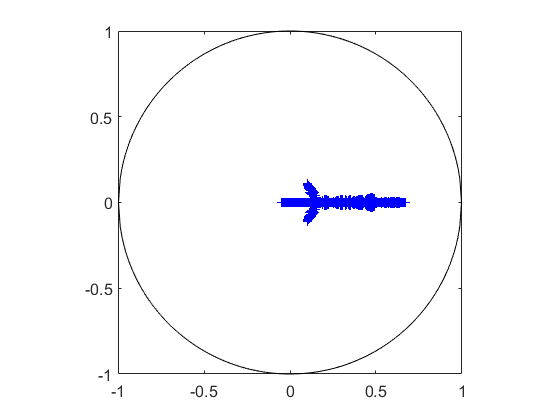}
\includegraphics[width=7.5cm,height=5.5cm]{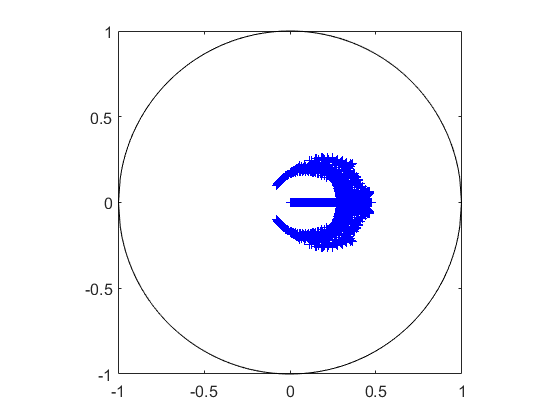}
\caption{The LFA-predicted spectra of the two-grid error-propagation operators for  the $P_2-P_1$ discretization. Left: VKI with $k=1$. Right: VKE with $k=1$.} \label{P2P1VK-TG}
\end{figure}

\begin{table}
\centering
\caption{Two-grid LFA predictions for Chebyshev  Vanka relaxation and two-grid performance with periodic ($\hat{\rho}$)  and Dirichlet ($\hat{\rho}^{\rm D}$) boundary conditions for the $P_2-P_1$ discretization for VKE. Optimal intervals for different $k$ with no weights. }
\begin{tabular}{l|l|l||l|l||l|l||l|l}
\hline
\multicolumn{3} {c||} {LFA predictions}  &\multicolumn{2} {c||} {$h=\frac{1}{20}$}   &\multicolumn{2} {c||} {$h=\frac{1}{40}$}  & \multicolumn{2} {c} {$h=\frac{1}{80}$} \\
\hline
   $k$  &$[\alpha,\beta]_{\rm VKE}$  &$\rho_{\rm VKE}$  & $\hat{\rho}_{\rm VKE}$  & $\hat{\rho}^{\rm D}_{\rm VKE}$   &$\hat{\rho}_{\rm VKE}$  & $\hat{\rho}^{\rm D}_{\rm VKE}$   &$\hat{\rho}_{\rm VKE}$  & $\hat{\rho}^{\rm D}_{\rm VKE}$ \\
\hline
 1                     &$[0.3,6.0]$              &0.475     &0.476   &0.571      &0.491  &0.571     &0.475  &0.571     \\
\hline
 2                     &$[0.5,4.7]$             &0.440     &0.412   &0.392      &0.485  &0.419       &0.388  &0.463     \\
\hline
 3                     &$[1.2,4.6]$             &0.168    &0.169    &0.175       &0.169   &0.176     &0.168  &0.176     \\
\hline
 4                     &$[1.6,4.6]$             &0.127    &0.125   &0.124        &0.128   &0.127     &0.127  & 0.127     \\
\hline
 5                     &$[2.6,3.7]$             &0.112    &0.111   &0.108        &0.111  &0.111      &0.112   &0.111    \\
\hline
\hline
\end{tabular}\label{No-TG-VKE}
\end{table}
\begin{remark}
In the numerical tests, we see that for some cases, the optimal pair of $[\alpha,\beta]$  are not unique. We break such ties arbitrarily.
\end{remark}
In Figure \ref{P2P1VK-SMOOTHING}, we show the LFA amplification factors using  Chebyshev Vanka relaxation with $k=1$ and the parameters from Tables \ref{No-TG-VKI} and \ref{No-TG-VKE}  in \eqref{Cheby-form-TG} for the $P_2-P_1$ discretization. We see that the VKE reduces the high-frequency error faster than VKI, and both reduce the low-frequency error slowly. In Figure \ref{P2P1VK-TG}, we present the spectrum of the associated two-grid error-propagation operators. The distribution of the eigenvalues is notably different: For VKI, most of the eigenvalues are real, while, for VKE, the eigenvalues are clustered around a circle in the complex plane.

\subsubsection{Geometric weights}\label{Natural-VK}

As seen in Tables \ref{No-TG-VKI} and  \ref{No-TG-VKE}, the measured convergence factors are largely independent of meshsize. Thus, we only consider $h=\frac{1}{40}$ in the following tests. Table \ref{W-MG} gives the measured convergence versus LFA predictions for relaxation with natural weights and  the corresponding optimal intervals. We see a good agreement between the LFA predictions and the measured convergence factors.  For the case of Dirichlet boundary conditions, the LFA predictions match the measured convergence factors very well except  for VKEW with $k=2$. Comparing Tables \ref{No-TG-VKI} and \ref{No-TG-VKE} with Table \ref{W-MG}, the big difference is that, for all $k$, VKIW outperforms VKEW. Moreover, for the inclusive patch, VKIW performs better than VKI. However, this is not true for the exclusive patch, where we see VKE has  better convergence factors than VKEW.

 \begin{table}
\centering
\caption{Two-grid LFA predictions for Chebyshev Vanka relaxation and  multigrid performance with periodic ($\hat{\rho}$) and Dirichlet ($\hat{\rho}^{\rm D}$) boundary conditions for the $P_2-P_1$ discretization. Optimal intervals for different $k$ with natural weights. $h=\frac{1}{40}$.}
\begin{tabular}{l|l|l|l|l||l|l|l|l}
\hline
\hline
 k       &$[\alpha,\beta]_{\rm VKIW}$     &$\rho_{\rm VKIW}$  & $\hat{\rho}_{\rm VKIW}$  & $\hat{\rho}^{\rm D}_{\rm VKIW}$  &$[\alpha,\beta]_{\rm VKEW}$  &$\rho_{\rm VKEW}$ &$\hat{\rho}_{\rm VKEW}$   &$\hat{\rho}^{\rm D}_{\rm VKEW}$ \\
\hline
 1                     &$[0.9, 2.9]$                      &0.518    &0.518  &0.517                 &$[1.3, 4.0]$            &0.584    &0.584   &0.589    \\
\hline
 2                     &$[1.1, 1.7]$                      &0.196    &0.197  &0.196                &$[0.5, 3.5]$            &0.376    &0.426    &0.279    \\
\hline
 3                     &$[1.4, 2.0]$                      &0.106    &0.126  &0.103                 &$[1.3, 3.6]$            &0.233    &0.234  &0.232      \\
\hline
 4                     &$[1.8, 2.2]$                      &0.085    &0.085   &0.085              &$[2.0, 3.5]$            &0.149     &0.149   &0.148     \\
\hline
 5                     &$[1.3,1.8]$                      &0.070    &0.071   &0.069              &$[2.5,3.5]$            &0.108     &0.107   &0.108      \\
\hline
\hline
\end{tabular}\label{W-MG}
\end{table}


\subsubsection{Optimized weights}
From the above results, we see that LFA provides a good prediction for the actual two-grid performance,  especially for the periodic problem. We also see that using geometric weights can improve performance. Motivated by this, we now  consider  whether using different weights for each different type of DoF within the relaxation scheme can improve performance. In this subsection, we apply LFA to optimize such  weights.  Here, we consider a preconditioned Richardson iteration,  with corresponding two-grid error propagation operator
\begin{equation}\label{Rich-error-form}
  E_R = (I-\omega_2 M_h^{-1}K_h)^{\nu_2}\mathcal{M}^{{\rm CGC}}(I-\omega_1 M_h^{-1}K_h)^{\nu_1}.
\end{equation}
Our target is to optimize the corresponding convergence factor, $\rho$, by brute-force search or using other optimization algorithms. We also consider the effect of using different pre- and post-relaxation weights in \eqref{Rich-error-form}.

Table \ref{Richardson-LFA} shows results for the $P_2-P_1$ discretization when optimizing only the outer parameters in the Richardson relaxation. We use brute-force search with sampling points taken in steps of 0.02 on the interval $[0,1]$ for the case of $\nu_1+\nu_2=1$. \revise{Fixing this weight, we then consider the performance of symmetric cycles with $\nu_1 = \nu_2 = 1,\ldots,5$, for comparison to the results using a Chebyshev iteration presented in Tables \ref{No-TG-VKI} through \ref{W-MG}.  We note both that the optimization of the Chebyshev intervals is more stable than simple use of multiple steps of a preconditioned Richardson iteration with fixed weight, and that the corresponding convergence factors presented in Table \ref{Richardson-LFA} are generally noticeably worse than those in Tables \ref{No-TG-VKI} through \ref{W-MG}, although perhaps further optimization of the Richardson weight would lead to some improvements.}  Further, for the case of $\nu_1+\nu_2=2$, we consider independent choices of $\omega_1$ and $\omega_2$, using sampling points taken in steps of 0.02 on the intervals $[0,0.5]$ and  $[0.5,0.9]$, respectively, to find the optimal results. Note that these intervals were selected based on results from a coarse sampling of  wider intervals.  We note that using $\nu_1+\nu_2=1$ seems to be more efficient than using $\nu_1+\nu_2=2$,  except for the case of VKI, even when using  two different weights for the pre- and post-relaxation parameters. We note that using different pre- and post-relaxation parameters gives notable improvements for VKE and VKI, but at best marginal gains for VKIW and VKEW.

Next, we consider fixing $\omega_1=\omega_2=1$ in \eqref{Rich-error-form} and using three different weights for $D_i$:  $d_1$ for  $N$-type velocity DoFs, $d_2$ for $X$-, $Y$-, and $C$-type  velocity DoFs,  and $d_3$ for the pressure. We make use of a robust optimization algorithm designed for LFA optimization\cite{LFAoptAlg} to find the optimal parameters, rather than brute force searches.  For the optimal parameters, we truncate the results obtained by the robust optimization to two  digits, noting that the performance is not overly  sensitive to this truncation. We first optimise for a single relaxation sweep, that is, $\nu_1+\nu_2=1$. Table \ref{3-weights-P2P1} shows that, in this setting,  Vanka inclusive  achieves a convergence factor of 0.581, while Vanka exclusive has a better convergence factor, 0.456. Both of these are significantly better than the corresponding results for VKI and VKE from Table \ref{Richardson-LFA}, but only slightly better than VKIW and VKEW.  Then, we optimise for the case of $\nu_1+\nu_2=2$, showing that there is no significant improvement compared with a single relaxation. Comparing Tables \ref{Richardson-LFA} and \ref{3-weights-P2P1} suggests that optimizing weights in $D_i$  obtains  better performance than doing so for the outer weights, although we note that nonsymmetric weighting of VKE in Table \ref{Richardson-LFA} outperforms the symmetric results in Table \ref{3-weights-P2P1}.

In order to see whether using more weights for different types of DoFs  can improve the performance of Vanka relaxation, we consider fixing $\omega_1=\omega_2=1$ and  using four weights, $d_1, d_2, d_3, d_4$, for the $N$-, $X$-, $Y$-, and $C$-type  velocity DoFs, respectively, and one weight, $d_5$,  for the pressure in $D_i$.  We again use the robust optimization algorithm. \cite{LFAoptAlg} Similarly,  we first consider a single relaxation. The optimal weights and corresponding LFA predictions are presented in Table \ref{5-weights-P2P1}, showing that VKE achieves  better performance than VKI. Then, we  optimise with $\nu_1+\nu_2=2$.  We achieve a convergence factor of 0.415 for Vanka inclusive and 0.408 for Vanka exclusive, which are only slightly  better than results  using three weights. Table  \ref{5-weights-P2P1} suggests that a single relaxation is again more efficient for VKI and VKE, especially for the case of VKE. All in all,  comparing Table \ref{3-weights-P2P1} with Table  \ref{5-weights-P2P1} shows that using three weights for $D_i$ is enough to obtain near-optimal performance. It is not necessary to use five weights.



\begin{table}
\centering
\caption{Two-grid LFA predictions, $\rho^{(\nu_1,\nu_2)}$, for Richardson relaxation optimizing outer weights for the $P_2-P_1$ discretization.}
\begin{tabular}{l||c|c|c|c|c|c|c||c|c}
\hline
\hline
 Method              &$\omega_{\rm opt}$   & $\rho^{(1,0)}$& $\rho^{(1,1)}$  &$\rho^{(2,2)}$    &$\rho^{(3,3)}$   &$\rho^{(4,4)}$   &$\rho^{(5,5)}$  &$(\omega_{1,\rm opt},\omega_{2,\rm opt})$   &$\rho^{(1,1)}$    \\
\hline
VKI                    & 0.24      &0.819       &0.670        &0.509      &0.411    &0.334         &0.273
&(0.14, 0.50)       & 0.556\\
\hline
VKIW                   &0.78       &0.587      &1.061         &0.337      &0.386    &0.107         &0.106
&(0.16, 0.84)       &0.507 \\
\hline
VKE                    &0.36       &0.669      &0.497        &0.638      &0.298    &0.126         &0.149
&(0.22, 0.56)        &0.356\\
\hline
VKEW                   &0.68       &0.574&1.507        &0.591      &0.919    &0.292         &0.473
&(0.00, 0.68)           &0.574 \\

\hline
\hline
\end{tabular}\label{Richardson-LFA}
\end{table}


\begin{table}
\centering
\caption{Two-grid LFA predictions for Richardson relaxation for the $P_2-P_1$ discretization with three weights in $D_i$.}
\begin{tabular}{l||c|c|c|l}
\hline
\hline
 Method              &$d_{\rm 1, opt}$   & $d_{\rm 2,opt}$  &$d_{3,\rm opt}$   &$\rho$    \\
\hline
VKI ($\nu_1+\nu_2=1$)                  &0.19               & 0.22             &0.71          &0.581    \\
VKE ($\nu_1+\nu_2=1$)                  &0.54               &0.26              &0.68          &0.456    \\
\hline
VKI ($\nu_1+\nu_2=2$)                  &0.22           &0.29          &0.47          &0.436    \\
VKE($\nu_1+\nu_2=2$)                   &0.39           &0.29          &0.37          &0.406   \\

\hline
\hline
\end{tabular}\label{3-weights-P2P1}
\end{table}

\begin{table}
\centering
\caption{Two-grid LFA predictions for Richardson relaxation for the $P_2-P_1$ discretization with five weights in $D_i$.}
\begin{tabular}{l||c|c|c|c|c|l}
\hline
\hline
 Method              &$d_{\rm 1, opt}$   & $d_{\rm 2,opt}$  &$d_{3,\rm opt}$  &$d_{4,\rm opt}$ &$d_{5,\rm opt}$ &$\rho$    \\
\hline
VKI ($\nu_1+\nu_2=1$)                  &0.19          &0.21           &0.17          &0.35   &0.74       &0.571       \\
VKE ($\nu_1+\nu_2=1$)                  &0.49          &0.25          &0.24           &0.30   &0.68       &0.452  \\
\hline
VKI ($\nu_1+\nu_2=2$)                 &0.22            &0.28         &0.28           &0.30   &0.51      &0.415\\
VKE ($\nu_1+\nu_2=2$)                 &0.43           &0.33          &0.30           &0.32   &0.36       &0.408  \\

\hline
\hline
\end{tabular}\label{5-weights-P2P1}
\end{table}

\subsubsection{Sensitivity of convergence factors to relaxation weights}
\begin{figure}
\centering
\includegraphics[width=6.cm,height=5.5cm]{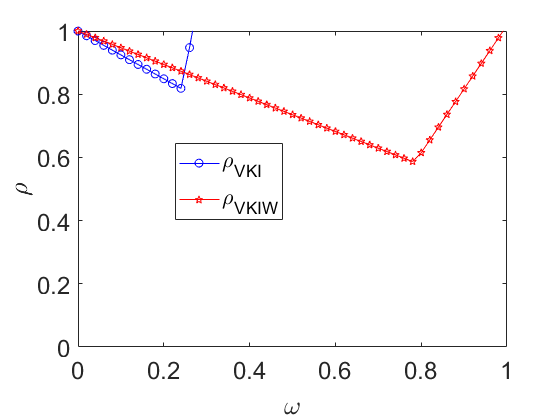}
\includegraphics[width=6.cm,height=5.5cm]{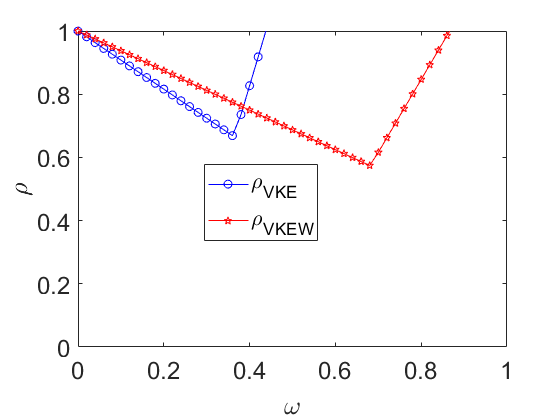}
\caption{Two-grid LFA convergence factor as a function of $\omega$ for Richardson relaxation with $\nu_1+\nu_2=1$. Left: Vanka-inclusive and Vanka-inclusive with natural weights. Right: Vanka-exclusive and Vanka-exclusive with natural weights.} \label{Rich-plot-sensitivity}
\end{figure}

Here, we present LFA results to show the sensitivity
of performance to parameter choice for Richardson relaxation with VKI, VKIW, VKE and VKEW. In Figure \ref{Rich-plot-sensitivity}, we show results for a single  Richardson relaxation for Vanka inclusive and Vanka exclusive, sampling  $\omega$ in steps of 0.02 on the interval $[0,1]$.  Figure \ref{Rich-plot-sensitivity} shows that it is better to underestimate the optimal parameter than to overestimate it. Similar behaviour is seen in other works.   Adler et al. \cite{adler2016monolithic} explored different types of Vanka as a preconditioner for GMRES for magnetohydrodynamics, showing the same preference for underestimation of the optimal parameter. In the application of Vanka relaxation for the Stokes equations with  an $H({\rm div})$ conforming discretization,\cite{JAdler_etal_2015b} even when using  different weights for the  velocity and pressure DoFs, it also appears better to  underestimate the optimal parameters.  We note that for both  Vanka-inclusive and Vanka-exclusive  with natural weights, convergence is observed for a wider range of Richardson weights than for  Vanka-inclusive and Vanka-exclusive.

\begin{figure}
\centering
\includegraphics[width=6.cm,height=5.5cm]{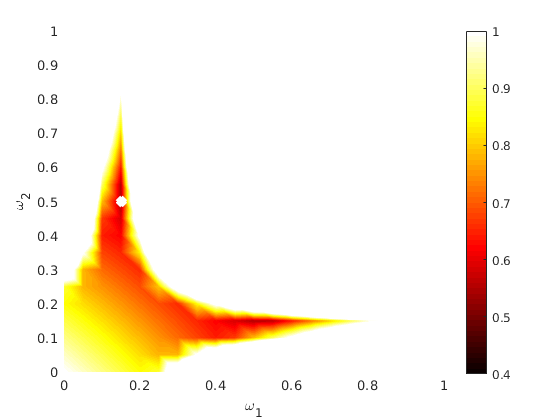}
\includegraphics[width=6.cm,height=5.5cm]{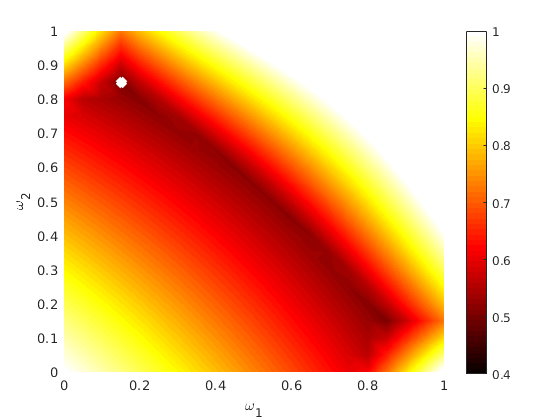}
\caption{Two-grid LFA convergence factor as a function of $\omega_1$ and $\omega_2$ for Richardson relaxation with $\nu_1+\nu_2=2$. The white circle marks the optimal point. Left: Vanka-inclusive. Right: Vanka-inclusive with natural weights.} \label{Rich-plot-2sweeps-VKI-sensitivity}
\end{figure}

\begin{figure}
\centering
\includegraphics[width=6.cm,height=5.5cm]{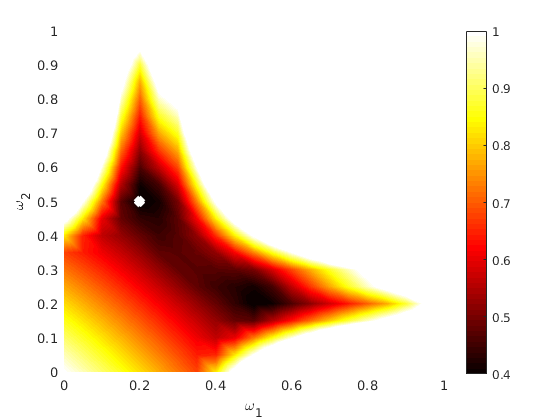}
\includegraphics[width=6.cm,height=5.5cm]{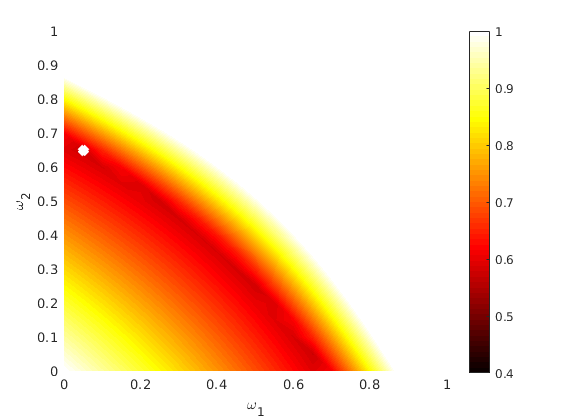}
\caption{Two-grid LFA convergence factor as a function of $\omega_1$ and $\omega_2$ for Richardson relaxation with $\nu_1+\nu_2=2$. The white circle marks the optimal points. Left: Vanka-exclusive. Right: Vanka-exclusive with natural weights.} \label{Rich-plot-2sweeps-VKE-sensitivity}
\end{figure}

In Figures \ref{Rich-plot-2sweeps-VKI-sensitivity} and \ref{Rich-plot-2sweeps-VKE-sensitivity}, we present the LFA convergence factor as function of $\omega_1$ and $\omega_2$, sampling in steps of 0.05 on the interval $[0,1]$, for Richardson relaxation with $\nu_1+\nu_2=2$.  We see similar behaviour, that performance degrades rapidly for weights that are ``too large''. For Vanka-inclusive, the case with no weights is more sensitive to the outer parameters, while the opposite is seen for  Vanka-exclusive.

\subsubsection{Sensitivity of convergence factors to mesh distortion}

\begin{figure}
\centering
\includegraphics[width=8.cm,height=5.5cm]{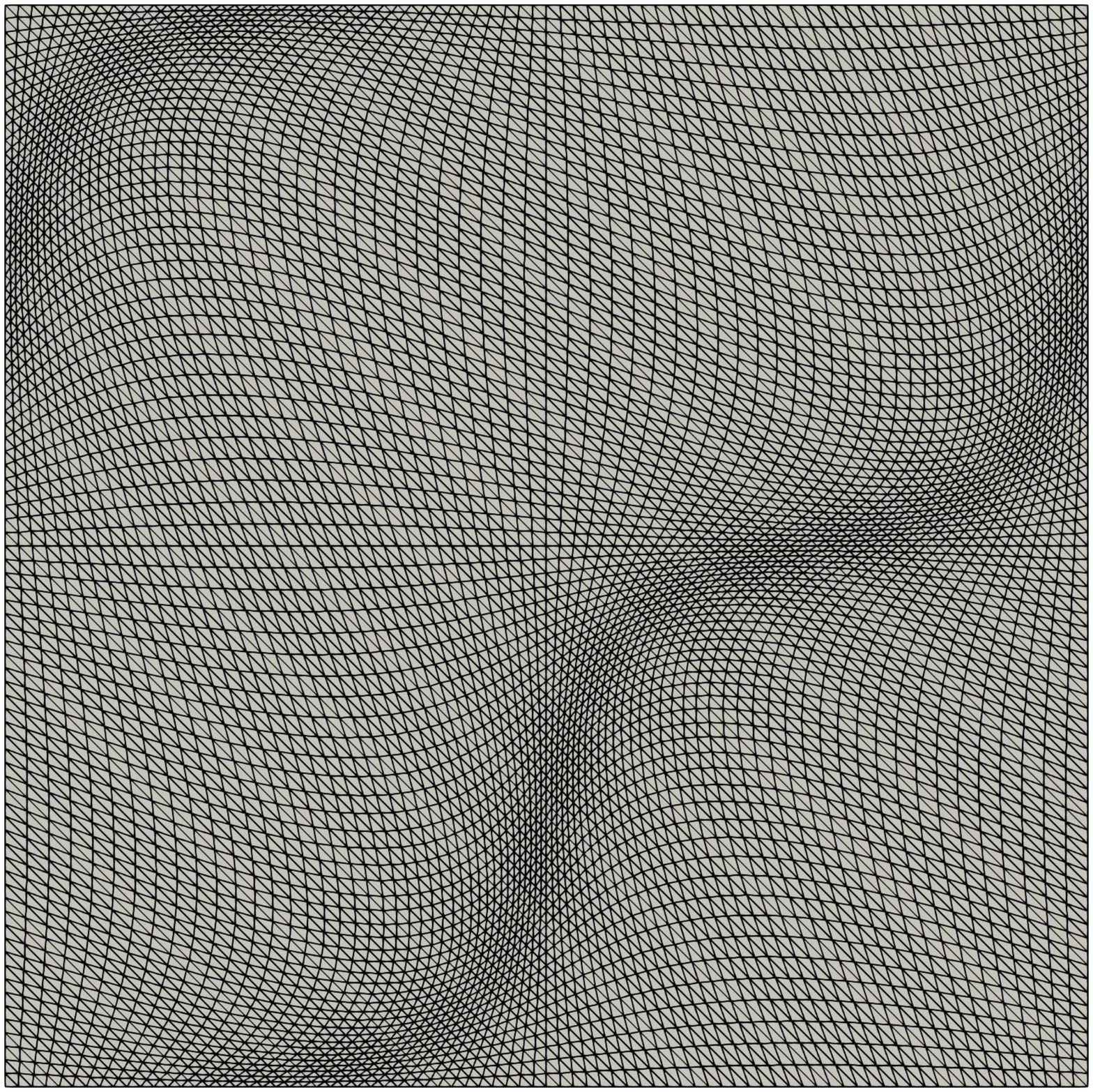}
\caption{Distorted mesh with $\varepsilon=0.1$ and $h=\frac{1}{80}$.} \label{distorted-mesh-eps0.1}
\end{figure}

\revise{Since local Fourier analysis applies only in the case of infinite or periodic meshes, a natural question is how accurate the predictions are in the case of distorted or otherwise non-uniform meshes.  Here, we consider smoothly distorted meshes, where the coordinates of the nodes of the (uniform) mesh are mapped under the transformation $(x,y) \rightarrow (x + \varepsilon\sin(2\pi x)\sin(2\pi y),y - \varepsilon\sin(2\pi x)\sin(2\pi y))$, and the connectivity of the mesh is preserved.  Figure \ref{distorted-mesh-eps0.1} shows the mesh distortion for the case of $\varepsilon=0.1$.  Tables \ref{P2P1-VKIE-distorted} and \ref{P2P1-VKIEW-distorted} show measured convergence factors for the two-grid algorithm as we distort a uniform mesh with $h=1/80$, taking $\varepsilon$ to be integer multiples of the meshsize, and fixing the Chebyshev relaxation parameters to be the same as the optimal choices for the uniform mesh.

Considering the results in Table \ref{P2P1-VKIE-distorted}, we see generally small changes in the measured convergence factor from those predicted by LFA for $\varepsilon = h = 0.0125$, although some changes, particularly for $k=2,3$ are more significant.  Slightly larger changes are seen in the case of $\varepsilon = 2h = 0.025$; however, for most cases the LFA predictions from the uniform mesh are surprisingly accurate.  Only for $k=2$ do we see significant deviation, measuring a convergence factor of 0.405 for VKI, compared to a prediction of 0.295, and 0.663 for VKE, compared to 0.440.  For larger $\varepsilon$, the discrepancy grows more significant (as expected), and we start to see true cases of ``failure'', with the optimal parameters from the uniform mesh leading to divergent stationary iterations.  Similar results are observed in Table \ref{P2P1-VKIEW-distorted} for the case of using geometric weights, with the most sensitivity to $\varepsilon$ seen for VKEW.}

 \begin{table}
\centering
\caption{Measured two-grid performance for multigrid with Chebyshev-Vanka relaxation with no weights applied to the $P_2-P_1$ discretization on distorted meshes. An alternating convergence pattern was observed, so the reported convergence factor is averaged over the final few iterations before
convergence (denoted *).  Diverging iterations are denoted by ``div''.}
\begin{tabular}{l|l|l|l|l|l||l|l|l|l|l}
\hline\hline
& \multicolumn{5} {c||} {VKI}  &\multicolumn{5} {c} {VKE} \\
\hline
\diagbox[width=3em]{k}{$\varepsilon$ }    &0       & 0.0125    &0.025  & 0.05     &0.1            & 0     &0.0125 &0.025 & 0.05  &0.1 \\
\hline
 1                                        &0.672   &0.673    &0.675  &0.679     &0.718            &0.475  &0.471   &0.479   &0.506*    &0.557* \\
\hline
 2                                        &0.295   &0.347    &0.405  &0.525    &0.755             &0.440  &0.540   &0.663   &div    &div \\
\hline
 3                                        &0.120   &0.176    &0.150*  &0.166    &0.266           &0.168    &0.210  &0.265   &0.393    &0.668 \\
\hline
 4                                         &0.102   &0.103   &0.104   &0.109   &0.127            &0.127    &0.141  &0.161   &0.184    &0.405 \\
\hline
 5                                         &0.051   &0.062    &0.079   &0.159    &0.316           &0.112   &0.117   &0.124   &0.139    &0.255 \\
\hline
\hline
\end{tabular}\label{P2P1-VKIE-distorted}
\end{table}

\begin{table}
\centering
\caption{Measured two-grid performance for multigrid with Chebyshev-Vanka relaxation with geometric weights applied to the $P_2-P_1$ discretization on distorted meshes. An alternating convergence pattern was observed, so the reported convergence factor is averaged over the final few iterations before
convergence (denoted *).  Diverging iterations are denoted by ``div''.}
\begin{tabular}{l|l|l|l|l|l||l|l|l|l|l}
  \hline
  \hline
& \multicolumn{5} {c||} {VKIW}  &\multicolumn{5} {c} {VKEW} \\
\hline
 \diagbox[width=3em]{k}{$\varepsilon$ }      &0     &0.0125 &0.025 & 0.05  &0.1                             & 0  &0.0125  &0.025 & 0.05  &0.1 \\
\hline
 1                                         &0.518   &0.540  &0.571  &0.636    &0.751                       &0.584  &0.586  &0.588   &0.607    &0.647 \\
\hline
 2                                         &0.196   &0.268  &0.382  &0.701     &div                        &0.376  &0.644  &div   &div    &div \\
\hline
 3                                          &0.106  &0.113*  &0.126  &0.228    &0.614                      &0.233   &0.302  &0.408   &0.740    &div \\
\hline
 4                                         &0.085   &0.118   &0.140  &0.169    &0.357                      &0.149   &0.175  &0.251   &0.440    &div \\
\hline
 5                                         &0.070   &0.088   &0.113   &0.168    &0.548                      &0.108  &0.153  &0.185   &0.256    &0.668 \\
\hline
\hline
\end{tabular}\label{P2P1-VKIEW-distorted}
\end{table}

\subsubsection{Comparison with multiplicative Vanka}

\revise{
  A recognized downside of additive relaxation schemes, as considered here, is that their convergence suffers in comparison to equivalent multiplicative schemes.  To some extent, this degradation in performance is seen as acceptable, since the additive schemes offer easy parallelization and are more efficient per iteration (since they require only one global residual evaluation per iteration, and not local residual updates at every stage). Nonetheless, it is worthwhile to understand how much performance is sacrificed when moving from a multiplicative to an additive scheme.  Table \ref{tab:Multiplicative} provides measured two-grid convergence factors for the finest-grid considered above with $h=1/80$, for a parallelized variant of multigrid with multiplicative Vanka relaxation run on 8 cores.  Here, we use the standard ``hybrid'' version of the relaxation scheme, which is multiplicative on a core, but additive across cores.

  To optimize the relaxation parameters, we use a two-stage brute force approach.  In the first stage, a coarse sampling of the endpoints of the Chebyshev interval was taken, sweeping over the lower endpoint in steps of 0.5, for values $0.5 \leq \alpha \leq 9.5$, and, for each value of $\alpha$, the upper endpoint in steps of 0.5, for values $\alpha+0.5 \leq \beta \leq 10.0$.  From the optimal parameters identified by this sweep, $[\alpha_c,\beta_c]$ (with lowest measured convergence factor at the iteration where the relative norm of the residual reaches $10^{-20}$), a finer sampling is used, taking steps of 0.1 in both $\alpha$ and $\beta$, from $\alpha_c-0.5$ to $\beta_c+0.5$, again ensuring $\alpha < \beta$.  The parameters reported in Table \ref{tab:Multiplicative} are those identified as optimal from this search and the convergence factors, $\rho$, reported are those measured at the iteration where the relative norm of the residual reaches $10^{-20}$; in cases where the optimal parameters occurred as boundary cases at the lower/upper limits of the search range, the search range was extended by 0.5 to resolve the possible ambiguity.

  Comparing performance between the multiplicative results in Table \ref{tab:Multiplicative} and those in Tables \ref{No-TG-VKI} and \ref{No-TG-VKE}, we see that the multiplicative form of the relaxation yields substantially better convergence, particularly for $k=1$ or $2$.  On the one hand, this emphasizes the common outcome for multiplicative schemes over their additive counterparts.  On the other hand, with increasing trends towards processors that excel at fine-grained parallelism, the fact that the additive schemes offer at least comparable performance is encouraging.
  }

\begin{table}
\centering
\caption{Measured two-grid performance for multigrid with multiplicative Chebyshev-Vanka relaxation with no weights, using brute-force optimization to determine the Chebyshev interval, applied to the $P_2-P_1$ discretization.}
\begin{tabular}{c||c|c||c|c||c|c||c|c}
  \hline
  \hline
  &  \multicolumn{2} {c||} {$k=1$}  & \multicolumn{2} {c||} {$k=2$}  & \multicolumn{2} {c||} {$k=3$}  & \multicolumn{2} {c} {$k=4$} \\
  & $[\alpha,\beta]$ & $\rho$   & $[\alpha,\beta]$ & $\rho$   & $[\alpha,\beta]$ & $\rho$   & $[\alpha,\beta]$ & $\rho$ \\
  \hline
  VKI & $[0.5,3.6]$ & 0.308 & $[1.2,3.0]$ & 0.068 & $[1.8,2.5]$ &  0.044 & $[2.2,2.3]$ & 0.031\\
  VKE & $[0.5,3.3]$ & 0.338 & $[1.1,2.7]$ & 0.075 & $[1.3,2.4]$ & 0.047 & $[1.5,2.1]$ & 0.038 \\
\hline
\hline
\end{tabular}\label{tab:Multiplicative}
\end{table}

\subsection{Numerical results for the $Q_2-Q_1$ discretization}
Similarly to the case of the $P_2-P_1$  discretization, we  consider Chebyshev-Vanka relaxation and preconditioned Richardson relaxation within monolithic multigrid methods for the Stokes equations for the $Q_2-Q_1$ discretization. Here, all of the optimal parameters are obtained by using the  robust optimization algorithm,\cite{LFAoptAlg} and we compare the LFA predictions with performance observed for the discretization with periodic boundary conditions  with $h=\frac{1}{40}$. We remark that the measured convergence factors with Dirichlet boundary conditions also match well with the LFA predictions.  For the optimal parameters, we truncate the results obtained by the robust optimization to two  digits, noting that the performance is not sensitive to this truncation.

Tables \ref{No-W-MG-Q2Q1} and \ref{W-MG-Q2Q1} present the measured convergence   for periodic boundary conditions versus  the LFA predictions for different $k$ and the corresponding optimal intervals for Chebyshev-Vanka relaxation with no weights and natural weights, respectively. Similarly to the case of the $P_2-P_1$ discretizations, we see that the optimal pair of $[\alpha,\beta]$ is not unique in some cases. We break such ties arbitrarily. We see  good agreement between the LFA predictions and the multigrid performance, except for VKI with $k=5$ and VKEW with $k=2$. For the former case, the per iteration convergence factor oscillates in the range from 0.12 to 0.22, averaging close to that predicted by LFA.  For the latter, more variation is seen in the per-cycle convergence factors, between 0.07 and 1.19, but the long-term  average still approaches that given by LFA.  For these cases, we present averaged convergence over the final five iterations. Comparing Tables \ref{No-W-MG-Q2Q1} and \ref{W-MG-Q2Q1} shows that the patches with natural weights outperform those with no weights, in some cases by a large margin.  An interesting observation from these tables, particularly Table \ref{W-MG-Q2Q1}, is that the optimization has a natural preference for giving very small intervals.  This would normally be a cause for concern, but clearly leads to excellent behaviour in this setting.

Table \ref{Richardson-LFA-Q2Q1} gives results for the $Q_2-Q_1$ discretization when optimizing only the outer parameter as in Table \ref{Richardson-LFA}.  \revise{As above, we then fix this weight and consider symmetric cycles with $\nu_1 = \nu_2 = 1,\ldots 5$.  For small numbers of sweeps, we see here that Richardson acceleration is as effective as Chebyshev for VKI and VKE, although Chebyshev is more effective for larger numbers of relaxation sweeps.  For VKIW and VKEW, in contrast, the Chebyshev-based relaxation schemes are always significantly more effective.} While using 2 different weights within a $(1,1)$ cycle is effective for VKI and VKE, better efficiency is seen for a single sweep with natural weights. We also consider  fixing $\omega_1=\omega_2=1$ in \eqref{Rich-error-form} and optimizing the inner weights. Table \ref{3-weights-Q2Q1} shows the optimal results for preconditioned Richardson relaxation with three parameters. Note that, as above, optimizing five weights does not improve the performance; thus, we omit the results here. Optimizing three weights obtains significantly better results than above for VKI, but shows little improvement when increasing the number of relaxation sweeps.  As above, optimizing three weights offers some improvement over natural weights, but not enough to suggest further optimization is worthwhile.

\begin{table}
\centering
\caption{Two-grid LFA predictions for Chebyshev Vanka relaxation and two-grid performance with periodic boundary conditions ($\hat{\rho}$)  for the $Q_2-Q_1$ discretization. Optimal intervals for different $k$ with no weights. $h=\frac{1}{40}$.  An alternating convergence pattern was observed, so the reported convergence factor is averaged over the final five iterations before
convergence (denoted *). }
\begin{tabular}{l|l|l|l|l||l|l}
\hline
\hline
 k  &$[\alpha,\beta]_{\rm VKI}$      &$\rho_{\rm VKI}$  & $\hat{\rho}_{\rm VKI}$    &$[\alpha,\beta]_{\rm VKE}$  &$\rho_{\rm VKE}$  &$\hat{\rho}_{\rm VKE}$  \\
\hline
 1                     &$[4.61, 4.81]$               &0.866    &0.864                &$[3.40, 3.49]$           &0.770     &0.769       \\
\hline
 2                      &$[4.72, 4.73]$              &0.752    &0.752                & $[0.96, 5.99]$         &0.506     &0.509       \\
\hline
 3                      &$[3.66, 5.82]$              &0.642    &0.635                & $[0.79, 6.18]$         &0.326     &0.327       \\
\hline
 4                      &$[0.39, 9.14]$              &0.203    &0.201                & $[0.91, 6.03]$           &0.262     &0.263       \\
\hline
 5                      &$[0.37, 11.42]$             &0.162    &0.166*             & $[3.29, 3.67]$           &0.278     &0.277       \\
\hline
\hline
\end{tabular}\label{No-W-MG-Q2Q1}
\end{table}

\begin{table}
\centering
\caption{Two-grid LFA predictions for Chebyshev Vanka relaxation and two-grid performance with periodic boundary conditions ($\hat{\rho}$)  for the $Q_2-Q_1$ discretization. Optimal intervals for different $k$ with natural weights. $h=\frac{1}{40}$.  An alternating convergence pattern was observed, so the reported convergence factor is averaged over the final five iterations before
convergence (denoted *).}
\begin{tabular}{l|l|l|l|l||l|l}
\hline
\hline
 k  &$[\alpha,\beta]_{\rm VKIW}$      &$\rho_{\rm VKIW}$  & $\hat{\rho}_{\rm VKIW}$    &$[\alpha,\beta]_{\rm VKEW}$  &$\rho_{\rm VKEW}$  &$\hat{\rho}_{\rm VKEW}$  \\
\hline
 1                     &$[0.15, 2.95]$              &0.637    &0.636                & $[1.32, 3.51]$           &0.681     &0.683       \\
\hline
 2                      &$[0.71, 1.97]$              &0.288    &0.290                & $[1.50, 1.54]$           &0.271     &0.265*        \\
\hline
 3                      &$[0.90, 1.42]$              &0.153    &0.161                & $[1.93, 1.99]$           &0.234     &0.241       \\
\hline
 4                      &$[1.23, 1.24]$              &0.097    &0.098                & $[2.38, 2.39]$           &0.217     &0.221       \\
\hline
 5                      &$[1.17, 1.40]$              &0.071    &0.072                & $[1.57, 2.41]$           &0.129     &0.128       \\
\hline
\hline
\end{tabular}\label{W-MG-Q2Q1}
\end{table}

\begin{table}
\centering
\caption{Two-grid LFA predictions, $\rho^{(\nu_1,\nu_2)}$, for Richardson relaxation optimizing outer weights for the $Q_2-Q_1$ discretization.}
\begin{tabular}{l||c|c|c|c|c|c|c||c|c}
\hline
\hline
 Method              &$\omega_{\rm opt}$   & $\rho^{(1,0)}$& $\rho^{(1,1)}$  &$\rho^{(2,2)}$    &$\rho^{(3,3)}$   &$\rho^{(4,4)}$   &$\rho^{(5,5)}$   &$(\omega_{1,\rm opt},\omega_{2,\rm opt})$   &$\rho^{(1,1)}$    \\
\hline
VKI                    &0.21      &0.931 &0.867        &0.752     &0.651    &0.565   &0.490
&(0.13, 0.60)   &0.809 \\
\hline
VKIW                   &0.92      &0.712 &0.913       &0.574      &0.290    &0.155   &0.134
&(0.65, 0.65)   &0.637 \\
\hline
VKE                    &0.29      &0.878 &0.770       &0.602     &0.488    &0.403   &0.337
&(0.76, 0.17)   &0.639 \\
\hline
VKEW                   &0.73      &0.697 &1.519       &0.507     &0.296    &0.577   &0.305
&(0.41, 0.41)   &0.684 \\

\hline
\hline
\end{tabular}\label{Richardson-LFA-Q2Q1}
\end{table}



\begin{table}
\centering
\caption{Two-grid LFA predictions for Richardson relaxation for the $Q_2-Q_1$ discretization with three weights in $D_i$.}
\begin{tabular}{l||c|c|c|l}
\hline
\hline
 Method              &$d_{\rm 1, opt}$   & $d_{\rm 2,opt}$  &$d_{3,\rm opt}$   &$\rho$    \\
\hline
VKI ($\nu_1+\nu_2=1$)                  &0.10       &0.13             &1.00         &0.695    \\
VKE ($\nu_1+\nu_2=1$)                  &0.88       &0.20             &0.84         &0.648    \\
\hline
VKI ($\nu_1+\nu_2=2$)                  &0.15       &0.09             &0.79         &0.583    \\
VKE($\nu_1+\nu_2=2$)                   &0.18       &0.22             &0.48         &0.646    \\
\hline
\hline
\end{tabular}\label{3-weights-Q2Q1}
\end{table}

\subsection{Comparing cost and performance}

\subsubsection{Cost of relaxation for the $P_2-P_1$ discretization}

 All schemes have the same cost for computing the initial residual. The $P_2$ Laplacian,  $A$, contains four types of stencils, that is, a 9-point stencil for $N$-type points, and 5-point stencils for $X$-,$Y$- and $C$-type points, and  $B_x$ and $B_y$ have 10-point stencils. Away from boundaries, we can naturally associate each DoF with a node, to see that a mesh with $n$ nodes also has about $n$ of each type of edge DoF.  So, the cost of a single residual evaluation on a mesh with $n$ nodes is (roughly) that of $(9+3\cdot5)n\cdot 2+10\cdot n\cdot4=88n$ multiply-add operations, coming from the 6 nonzero blocks in the matrix.

For the $P_2-P_1$ discretization,  the remaining  cost is that of  solving a small  problem, $K_ix_i=b_i$, in each patch. Here, we use $LU$ decomposition to solve these subproblems. Assume that the $LU$ decomposition is precomputed, as this can be done once per patch and  used for all of the solves over that patch. As in Table \ref{Patch-points-number},  for VKI and VKIW,  the size of the patch problem is $2(7+4+4+4)+1=39$.  For VKE and VKEW, it is $2(1+2+2+4)+1=27$. Note that, similarly to $K$, $K_i$ is a block (sparse) matrix, and its $LU$ factors will retain some sparsity as well.  Thus, the cost of solving the sparse systems with $L$ or $U$  will  require a number of multiply-add operations equal to the number of nonzero entries in $L$ or $U$.  Table \ref{Patch-DoFs-number} presents these numbers, based on direct calculation of the factorizations. For VKI and VKIW, the cost  of applying the inverses of $L$ and $U$ is $566$ multiply-add operations. For the whole system, we need to solve roughly  $n$ subproblems, giving a total cost of $566n$. Similarly, it costs $230n$ multiply-add operations for VKE and VKEW.  For the approaches with natural weights, we need an additional scaling for each subproblem. Thus, there will  be an additional cost of $39 n$ and $27n$ for VKIW and VKEW, respectively.

\begin{table}
\centering
\caption{The number (\#) of  nonzero elements of $L$ and $U$ ($LU=K_i$) for Vanka-inclusive and Vanka-exclusive patches for the $P_2-P_1$ and $Q_2-Q_1$ discretizations.}
\begin{tabular}{l||c|c||c|c}
   &\multicolumn{2} {c||} {$P_2-P_1$}  & \multicolumn{2} {c} {$Q_2-Q_1$} \\
\hline
 \backslashbox{\#}{patch}            &$\Xi_{VKI}$   & $\Xi_{VKE}$  &$\Xi_{VKI}$   &$\Xi_{VKE}$   \\
\hline
$L$                 &284     &115        &617      &335   \\
\hline
$U$                 &282     &115        &617       &337  \\
\hline
Total             &566      &230        &1234       &672 \\
\hline

\hline
\end{tabular}\label{Patch-DoFs-number}
\end{table}

Accumulating the costs of a residual evaluation with these, we have total costs of $88n+566n=654n$ multiply-add operations per sweep of VKI, $88n+230n=318n$ multiply-add operations per sweep of VKE,  $654n+39n=693n$ multiply-add operations per sweep of VKIW, and $318n+27n=345n$ multiply-add operations per sweep of VKEW. To compare these costs, we omit the cost of the coarse-grid correction and only consider the cost of the relaxation scheme. We denote the above cost as $W_s$, which corresponds to $k=1$ and $\nu_1+\nu_2=1$. Comparing efficiencies can now be easily done by appropriately
weighting either measured or predicted convergence factors relative to their work (here, we use the predicted convergence factors): if one iteration costs $W$
times that of another, and yields a convergence factor of $\rho_1$, then we can easily compare $\rho_1^{1/W}$
directly to the second convergence factor, $\rho_2$, to see if the effective error reduction achieved by
the first algorithm in an equal amount of work to the second is better or worse than that achieved
by the second. Here, we compare the efficiency relative to VKEW.

 Next,  using the data presented above, we  find the most effective parameters for each  of VKI, VKE, VKIW and VKEW, and compare their costs. From Tables \ref{No-TG-VKI} and \ref{No-TG-VKE}, we see the most effective choice for VKI is $k=3$ giving $\rho=0.120$, while, for VKE, it is $k=1$ giving $\rho=0.475$. Note, however, that a more efficient set of parameters for VKE is found in Table \ref{Richardson-LFA}, giving a convergence factor of 0.356 per iteration. From Table \ref{W-MG}, the most effective choice for VKIW is $k=2$ giving $\rho=0.196$, and $k=1$ giving $\rho=0.584$ for VKEW. Here, optimizing with more weights, as in Tables \ref{3-weights-P2P1} and \ref{5-weights-P2P1}, was more effective, and we consider the best results with a non-trivial scaling matrix, $D_i$, in place of the geometric weights, giving 0.571 for VKIW and 0.452 for VKEW.
 Now, we can calculate the efficiency for each relaxation scheme relative to its total cost, $W_t=k\cdot (\nu_1+\nu_2)\cdot W_s$, where $W_s$ denotes the cost of a single relaxation.   Table \ref{P2P1-COST} details the  cost-effectiveness for each scheme for the $P_2-P_1$ discretization, showing that VKEW offers the most efficient relaxation scheme. In particular, this shows that the substantially smaller cost per iteration of the exclusive patches offers greater efficiency, despite the improved convergence when using inclusive patches.

\begin{table}
\centering
\caption{Comparing cost and performance for the $P_2-P_1$ discretization.}
\begin{tabular}{l||l|l||l|l}
\hline
\hline
                               &VKI        & VKE          &VKIW        &VKEW   \\
\hline
$W_s$ ($k=1,\nu_1+\nu_2=1$)     &654n      &318n          &693n         &345n \\
\hline
Most Effective $(k,\nu_1+\nu_2)$    &(3, 2)   &(1, 2)         & (1, 1)           & (1, 1) \\
\hline
$\rho$                    &0.120      &0.356       &0.571        &0.452 \\
\hline
    \multicolumn{4} {r}  {$W_t = k\cdot (\nu_1+\nu_2) \cdot W_s$}\\
\hline
$W_t$                       &3924n      &636n        &693n         &345n \\
\hline
Relative Efficiency         &0.830      &0.571      &0.757         &0.452 \\
\hline
\hline
\end{tabular}\label{P2P1-COST}
\end{table}

\subsubsection{Cost of relaxation for the $Q_2-Q_1$ discretization}
 The $Q_2$ Laplacian, $A$, contains four types of stencils, a 25-point stencil for $N$-type points, 15-point stencils for $X$-, and $Y$-  DoFs and a 9-point stencil for $C$-type points, while $B_x$ and $B_y$ have 12-point stencils. So, the cost of a single residual evaluation on a mesh with $n$ nodes is (roughly) that of $(25+2\cdot15+9)n\cdot 2+12\cdot n\cdot4=104n$ multiply-add operations, coming from the 6 nonzero blocks in the matrix.

For the $Q_2-Q_1$ approximation with  additive Vanka-type relaxation, there is  little difference in the cost calculation compared with that of the $P_2-P_1$ discretization.  Here, we again use the $LU$ decomposition of the patch matrices, $K_i$, and assume the $LU$ decomposition is done.    For VKI,  the cost  of applying the inverses of $L$ and $U$ is $1234$ multiply-add operations per block, see Table \ref{Patch-DoFs-number}. Thus, the cost for a full sweep is $1234n$ multiply-add operations. Similarly, it costs $672n$ multiply-add operations per sweep of VKE.  For the approaches with natural weights,  there will  be an additional cost of $51n$ and $35n$ for VKIW and VKEW, respectively, for these scaling operations. Accumulating the costs of a residual evaluation with these, we have total costs of $104n+1234n=1338n$ multiply-add operations per sweep of VKI, $104n+672n= 776n$ multiply-add operations per sweep of VKE,  $1338+51n= 1389n$ multiply-add operations per sweep of VKIW, and $776n+ 35n= 811n$ multiply-add operations per sweep of VKEW.

From Tables \ref{No-W-MG-Q2Q1} and  \ref{W-MG-Q2Q1}, we see the most effective choice for VKI is $k=4$ giving $\rho=0.203$, while, for VKE, it is $k=3$ giving $\rho=0.326$. Note, however, that a more efficient set of parameters for VKE is found in Table \ref{Richardson-LFA-Q2Q1}, giving a convergence factor of 0.639 per iteration with $\nu_1+\nu_2=2$. From Table \ref{W-MG-Q2Q1}, the most effective choice for VKIW is $k=3$ giving $\rho=0.153$, and $k=2$ giving $\rho=0.271$ for VKEW. Here, optimizing with more weights, as in Table \ref{3-weights-Q2Q1}, was more effective, and we consider the best results with a non-trivial scaling matrix, $D_i$, in place of the geometric weights, giving 0.695 for VKIW and 0.648 for VKEW with $\nu_1+\nu_2=1$. Table \ref{Q2Q1-COST} compares the efficiency relative to VKEW, showing that the most efficient choice is VKEW,  the same as for the $P_2-P_1$ discretization.

\begin{table}
\centering
\caption{Comparing cost and performance for the $Q_2-Q_1$ discretization.}
\begin{tabular}{l||l|l||l|l}
\hline
\hline
                               &VKI        & VKE          &VKIW        &VKEW   \\
\hline
$W_s$ ($k=1,\nu_1+\nu_2=1$)     &1338n      &776n         &1389n         &811n \\
\hline
Most Effective $(k,\nu_1+\nu_2)$    &(4, 2)   &(1, 2)     & (1, 1)      & (1, 1) \\
\hline
$\rho$                           &0.203      &0.639        &0.695         &0.648 \\
\hline
    \multicolumn{4} {r}  {$W_t = k\cdot (\nu_1+\nu_2) \cdot W_s$}\\
\hline
$W_t$                        &10704n      &1552n         & 1389n        &811n \\
\hline
Relative Efficiency          &0.886      &0.791         &0.809         &0.648 \\
\hline
\hline
\end{tabular}\label{Q2Q1-COST}
\end{table}

\section{Conclusions and future work}\label{sec:concl-FW}

We present a local Fourier analysis for a  monolithic multigrid method based on overlapping additive Vanka-type relaxation for the Stokes equations. Two choices of  patches for the overlapping schemes are discussed for the $P_2-P_1$ and $Q_2-Q_1$ discretizations. A general framework of LFA for additive Vanka relaxation is developed to help choose algorithmic parameters, which can be applied to other problems and to different discretizations.  The LFA shows that using smaller patches can outperform relaxation using bigger patches due to the lower cost per sweep of relaxation. Moreover, to improve the performance, we use LFA to optimize the weights, yielding notable  improvement. Numerical performance with periodic and Dirichlet boundary conditions validate the LFA predictions, showing  that these Vanka relaxation schemes are robust to the different boundary conditions.

Extending Vanka relaxation for other types of problems is an interesting topic. We note that this LFA framework of additive Vanka relaxation has a limitation: since we need to know the patch first, selection of the patches is not readily optimized in this framework. Another interesting question is the use of other boundary conditions on the patches, such as are used in optimized Schwarz, which could also be tuned using LFA. Developing general-purpose LFA software for additive Vanka with automatic evaluation of patch choices and optimizing the weights are also topics for future work.
\section*{Acknowledgments}

The work of P.~E.~F.~was supported by the Engineering and Physical Sciences Research
Council [grant numbers EP/K030930/1 and EP/R029423/1]. The work of S.~P.~M.~was partially supported by an NSERC Discovery Grant.

\bibliography{stokes_refs}%

\end{document}